\documentclass[preprint,12pt]{elsarticle}

\usepackage{amssymb}

\usepackage[utf8]{inputenc}
\usepackage[T1]{fontenc}
\usepackage{lineno}
\modulolinenumbers[5]
\usepackage[english]{babel}
\usepackage{geometry}
\usepackage{multirow}
\usepackage[table,xcdraw]{xcolor}
\usepackage{times}
\usepackage{graphicx}
\usepackage{amsfonts}
\usepackage{latexsym}
\usepackage{amsmath}
\usepackage{mathtools}
\usepackage{array}
\usepackage{cellspace}
\usepackage{multicol}
\usepackage{enumerate}
\usepackage{setspace}
\usepackage{wrapfig} 
\usepackage{subcaption}

\newcommand{\ud}{\,\mathrm{d}}

\makeatletter
\def\ps@pprintTitle{%
    \let\@oddhead\@empty
    \let\@evenhead\@empty
    \def\@oddfoot{\footnotesize\itshape
         {A preprint} \hfill {March 4, 2022}}%
    \let\@evenfoot\@oddfoot
    }
\makeatother

\begin{document}

\begin{frontmatter}

\title{High-Degree Splines from Discrete Fourier Transforms: Robust Methods to Obtain the Boundary Conditions}

%
%
%
            
\author{A. Pepin$^1$, S. Léger$^2$, N. Beaudoin$^3$}

\address{$^1$Secteur sciences,
         Universit{\'e} de Moncton, Campus de Shippagan\\
         218 Bd. J. D. Gauthier,
         Shippagan, Canada,
         E8S 1P6.\\ }

\address{$^2$D{\'e}partement de math{\'e}matiques
         et de statistique,
         Pavillon Rémi-Rossignol,\\
         18 avenue Antonine-Maillet,
         Universit{\'e} de Moncton,
         Moncton, Canada,
         E1A 3E9.\\ }
         
\address{$^3$D{\'e}partement de physique et d'astronomie,
         Pavillon Rémi-Rossignol,\\
         18 avenue Antonine-Maillet,
         Universit{\'e} de Moncton,
         Moncton, Canada,
         E1A 3E9.\\ }

\begin{abstract}
Computing accurate splines of degree greater than three is still a challenging task in today's applications. In this type of interpolation, high-order derivatives are needed on the given mesh. As these derivatives are rarely known and are often not easy to approximate accurately, high-degree splines are difficult to obtain using standard approaches.

In Beaudoin (1998), Beaudoin and Beauchemin (2003), and Pepin \textit{et al.} (2019), a new method to compute spline approximations of low or high degree from equidistant interpolation nodes based on the discrete Fourier transform is analyzed. The accuracy of this method greatly depends on the accuracy of the boundary conditions. An algorithm for the computation of the boundary conditions can be found in Beaudoin (1998), and Beaudoin and Beauchemin (2003). However, this algorithm lacks robustness since the approximation of the boundary conditions is strongly dependant on the choice of $\theta$ arbitrary parameters, $\theta$ being the degree of the spline.

The goal of this paper is therefore to propose two new robust algorithms, independent of arbitrary parameters, for the computation of the boundary conditions in order to obtain accurate splines of any degree. Numerical results will be presented to show the efficiency of these new approaches.  
\end{abstract}


\end{frontmatter}


\section{Introduction}
In today's modern applications, the numerical methods developed need to be accurate, reliable and robust. This is especially true in the field of interpolation, as certain applications, such as signal processing, require very high accuracy \cite{unser_splines:_1999}. As high-degree polynomial interpolation usually leads to important oscillations at the edges of the interval, a phenomenon known as Runge's phenomenon, spline interpolation is typically preferred. 
Widely used in industry and research, spline interpolation has proved to be a very efficient tool, particularly when dealing with data interpolation or curve smoothing applications. In practice, cubic spline interpolation is commonly used due to its good accuracy to computational cost ratio. Many authors have studied the behaviour and theoretical context of this powerful tool \cite{ahlberg_theory_1967,de_boor_practical_1978,holladay_smoothest_1957,schoenberg_contributions_1946}. In particular, \cite{ahlberg_theory_1967} gives a more general approach for higher odd-degree spline interpolation. High-degree spline interpolation is however rarely used in practice as some authors believe it fails to yield any advantage (see \cite{cheney_numerical_2012} for example).
However, as shown in \cite{beaudoin_tutorial/article_1998}, computing an accurate numerical Fourier transform requires very accurate interpolation results, which means that efficient algorithms that lead to accurate high-degree splines are in fact needed.

It is known that even-degree splines behave differently than odd-degree splines \cite{ahlberg_theory_1967,cheney_numerical_2012}. As even-degree splines may not exist for a given mesh \cite{ahlberg_theory_1967}, are less efficient when interpolating at the nodes \cite{cheney_numerical_2012} and can lead to an inaccurate approximation of the digitized function (for example, when using a quadratic spline, as shown in \cite{de_boor_practical_1978}), odd-degree splines are typically preferred. Computing high-degree splines, no matter if they are of even or odd degree, does however come with challenges. The biggest challenge remains the computation of accurate boundary conditions without the need of additional information at the nodes.

In \cite{beaudoin_tutorial/article_1998, beaudoin_new_2003}, Beaudoin and Beauchemin established a new method to construct splines of any odd degree. Starting from a discrete function, an accurate approximation of the continuous Fourier transform was developed. It was then later formally demonstrated that the piecewise continuous polynomials defined in their work form a spline function and that this was the case for both odd and even degree polynomials \cite{pepin_new_2019}.

Several authors \cite{froeyen_improved_1985,makinen_new_1982,schutte_new_1981,sorella_improved_1984} have stated that the discrete Fourier transform (DFT) is far from being adequate to act in place of the Fourier transform. The approach proposed by Beaudoin and Beauchemin in \cite{beaudoin_tutorial/article_1998,beaudoin_new_2003} yields better results than the classical DFT, while providing additional information, such as accurate numerical derivatives and integration of the digitized function. With these numerical derivatives, it is then easy to obtain accurate high-degree splines, something that is generally difficult to achieve using a standard approach. In \cite{beaudoin_tutorial/article_1998}, Beaudoin and Beauchemin noted that the numerical accuracy of their method increases much more rapidly than the computational cost. However, just like other methods that construct spline functions, the accuracy of the results depends on the accuracy of the boundary conditions. Obtaining accurate approximations of the boundary conditions is therefore key to the good performance of the spline interpolation method. 

Perhaps its most interesting advantage, splines of degree greater than three can yield accurate approximations of the higher order derivatives. Such approximations are particularly useful when dealing with boundary value problems (see \cite{wazwaz_approximate_2000}), which has been the object of many studies throughout the years. Amongst others, Usmani has studied the use of quartic splines for the solution of a fourth-order boundary value problem in \cite{usmani_use_1992}, while Siddiqi and Twizell have studied the use of splines of degree six, eight, ten and twelve for the solution of linear boundary value problems of degree six, eight, ten and twelve  respectively \cite{siddiqi_spline_1996,siddiqi_spline_1996-1,siddiqi_spline_1998,siddiqi_spline_1997}.

In our proposed approach, as was already stated, the computation of the boundary conditions is a key ingredient. This paper therefore presents two new methods for the computation of the boundary conditions of splines of any degree. The method established by Beaudoin and Beauchemin, combined with our algorithms to compute the boundary conditions, will then be compared with the traditional cubic spline (for both natural and not-a-knot boundary conditions) to show its numerical efficiency.

The paper is organized as follows. Section~\ref{sec:MathBack} gives a brief summary of Beaudoin and Beauchemin's method. Section \ref{sec:AlgBC} presents our proposed algorithms for the approximation of the boundary conditions, while Section \ref{Sec:Numerical_results} presents some numerical results to show the good performance, in terms of accuracy and robustness, of our proposed algorithms. This paper will consider one-dimensional functions and, although these methods can be used for splines
of degree as high as one wants, will focus on spline
interpolation of degree 3, 4 and 5 for comparison purposes
with methods commonly used in practice. An example
involving splines of degree 7, 9 and 11 is presented in Section \ref{Sec:Practical_example}. The gain in accuracy for higher degree splines will also be shown. All numerical results will be obtained using MATLAB\textsuperscript{\textregistered} and Maple.  

\section{Mathematical background} \label{sec:MathBack}
In the following, $\displaystyle\mathbb{C}^\theta_\text{c}$ denotes the set of complex-valued column vectors of dimension $\theta$. Similarly, $\displaystyle\mathbb{C}^\theta_\text{r}$ denotes the set of complex-valued row vectors of dimension $\theta$.

Let us consider a real-valued analytic function $g$ defined on an interval $[0,T]$ for which its values are only known at $N+1$ equally spaced nodes $t_j=j\Delta t$, $j=0,1,\ldots,N$, with $\displaystyle\Delta t=T/N$. The value of $g$ at the node $t_j$ is denoted by $g(t_j)=g_j$. 

Starting from $N+1$ points, the goal is to get an approximation of the function $g$ by spline interpolation of degree $\theta\in\mathbb{N}$. The resulting spline function will be denoted by $g_\theta$.

A new method to compute odd-degree splines was developed by Beaudoin and Beauchemin in \cite{beaudoin_tutorial/article_1998,beaudoin_new_2003}. Their results led them to believe that this method was not applicable for even-degree splines. It was however
later shown in \cite{pepin_new_2019} that this method works for splines of any degree, as long as $N$ and $\theta$ are not simultaneously chosen as even integers.

The spline functions are obtained by defining a polynomial of degree $\theta$, denoted by $\left[g_\theta\right]_j$, on every subinterval $[t_j,t_{j+1}]$. Using the truncated Taylor development of the analytic function $g$ at the interpolation nodes $t_j$, the polynomial $\left[g_\theta\right]_j:[t_j,t_{j+1}]\rightarrow\mathbb{R}$ of degree $\theta$ on the subinterval $[t_j,t_{j+1}]$ is defined by
\begin{equation}
	\left[g_\theta\right]_j(t) = \sum_{\mu=0}^{\theta} \frac{(t-t_j)^\mu }{\mu!}g_{j,\theta}^{(\mu)}, \qquad\text{for}~t\in [t_j,t_{j+1}],
    \label{Polynome_fct_g}
\end{equation}
where $g_{j,\theta}^{(\mu)}$ is the $\mu$-th numerical derivative of order $\theta$ of $g$ evaluated at the node $t=t_j$. For $\mu=0$, the values $g_{j,\theta}$ are simply the values at the interpolation nodes (i.e. $g_{j}$).

From Equation \eqref{Polynome_fct_g}, the derivatives of $g_\theta$ are then obtained using the following equation
\begin{equation*}
	\left[g^{(\beta)}_\theta\right]_j(t) = \sum_{\mu=0}^{\theta-\beta} \frac{(t-t_j)^\mu}{\mu!}g_{j,\theta}^{(\mu+\beta)}, \qquad\text{for}~t\in [t_j,t_{j+1}],
\end{equation*}
where $\beta=1,2,\ldots,\theta-1$. Evidently, $g_{\theta}^{(\theta)}$ is not continuous on $[0,T]$ since it is defined by the constant $g_{j,\theta}^{(\theta)}$ on the interval $[t_j,t_{j+1})$.

In order to compute the numerical derivatives $g_{j,\theta}^{(\mu)}$, $\mu=1,2,\ldots,\theta$, for any positive integer $\theta$ using Beaudoin and Beauchemin's method, the following system must be solved for $F_{\theta,k}$
\begin{equation} 
	M_{\theta,k} F_{\theta,k} = B+C_k,
    \label{Eq_SystMatricielFb}
\end{equation}
where $M_{\theta,k}$ is a square matrix of dimension $\theta$ and $\displaystyle F_{\theta,k}$, $\displaystyle B$ and $\displaystyle C_k$ are column vectors of dimension $\theta$.

\noindent In Equation \eqref{Eq_SystMatricielFb}, $\displaystyle F_{\theta,k}\in\mathbb{C}^\theta_\text{c}$ is simply the column vector containing the discrete Fourier transform of the sequence $g_{j,\theta}^{(1)},g_{j,\theta}^{(2)},\ldots,g_{j,\theta}^{(\theta)}$ for $j=0,1,\ldots,N-1$ and is defined by
\begin{equation*}
	(F_{\theta,k})_{\mu,1} = f_{\mu,k,\theta} = \sum_{j=0}^{N-1}g_{j,\theta}^{(\mu)}\exp(-\mathrm{i}2\pi kj/N), \qquad\text{for~}\mu=1, 2, \ldots, \theta.
\end{equation*}

\noindent As for $M_{\theta,k}$, it is defined by
\begin{equation*}
	(M_{\theta,k})_{\mu,\nu}=\begin{cases}
      0, & \text{if}\ \nu-\mu+1<0,\\
      \displaystyle J_{\nu-\mu+1,k}, & \text{otherwise},\\
    \end{cases}
\end{equation*}
for $\mu,\nu=1,2,\ldots,\theta$, where
\begin{equation*}
	\displaystyle J_{p,k}=
    \begin{cases}
      \displaystyle \exp(-\mathrm{i}2\pi k/N)-1, & \text{if}\ p= 0,\\
      \displaystyle\frac{(\Delta t)^{p}}{p!}\exp(-\mathrm{i}2\pi k/N), & \text{if}\ p >0.\\
    \end{cases}
\end{equation*}
From \cite{pepin_new_2019}, the determinant of $M_{\theta,k}$ is given by
\begin{equation}
       	\det(M_{\theta,k})=\displaystyle\frac{(\Delta t)^{\theta}}{\theta!}\sum_{q=0}^{\theta-1} A_{\theta,q}\left(\exp(-\mathrm{i}2\pi k/N)\right)^{\theta-q},
        \label{Eq_determinant_Mb}
    \end{equation}
where $\displaystyle A_{\theta,q}$, for $0\leq q\leq\theta-1$, is the eulerian number defined by
\begin{equation*}
    A_{\theta,q}=\sum_{k=0}^{q}(-1)^k\frac{(\theta+1)!}{(\theta+1-k)!k!}(q+1-k)^\theta.
\end{equation*}
For $\theta=0$, we have $\det(M_{0,k})=1$. 

The proof for Equation \eqref{Eq_determinant_Mb} is given in \cite{pepin_new_2019}. It is also shown that Equation \eqref{Eq_determinant_Mb} is equal to zero if and only if $\theta$ is even and $k=N/2$. To avoid singularities in $M_{\theta,k}$ when $k=N/2$, it is clear that $\theta$ and $N$ must not be chosen simultaneously as even integers.

\noindent Furthermore, since $M_{\theta,k}$ is an upper Toeplitz-Hessenberg matrix, the explicit inversion formula for these kind of matrices (as shown in \cite{marrero_explicit_2018}) can straightforwardly be used to demonstrate that the elements of $M_{\theta,k}^{-1}$ are given by the following equation
\begin{equation}
	\left(M_{\theta,k}^{-1}\right)_{\mu,\nu} = 
    \begin{cases}
		\displaystyle(-J_{0,k})^{\mu-\nu}\left(\frac{\det(M_{\nu-1,k})\det(M_{\theta-\mu,k})}{\det(M_{\theta,k})}-\det(M_{\nu-\mu-1,k})\right), & \text{if}\ \nu>\mu,\\ 
        \displaystyle(-J_{0,k})^{\mu-\nu}\frac{\det(M_{\nu-1,k})\det(M_{\theta-\mu,k})}{\det(M_{\theta,k})}, & \text{if}\ \nu\leq\mu,
	\end{cases}
	\label{Eq_MatrixInverseMtheta}
\end{equation}
for $\mu,\nu=1,2,\ldots,\theta$. This equation gives a simple and explicit relation to compute the elements of $M_{\theta,k}^{-1}$ for any given value $\theta\in\mathbb{N}$. However, Equation \eqref{Eq_MatrixInverseMtheta} leads to an indeterminate form when computing the elements located above the main diagonal when $k=0$ (since $k=0$ implies $J_{0,k}=0$). 

While numerical observations have lead us to conclude that these indeterminate forms disappear upon simplification of Equation \eqref{Eq_MatrixInverseMtheta}, it is not easy to obtain a more general equation than the one presented. For this reason and without loss of generality, we will limit the use of Equation \eqref{Eq_MatrixInverseMtheta} to the case $\nu\leq\mu$.

The vector $\displaystyle B\in\mathbb{R}^\theta_\text{c}$ in Equation \eqref{Eq_SystMatricielFb}, which is typically unknown for a given digitized function $g$, is called the boundary conditions vector and is defined by
\begin{equation*}
	(B)_{\mu,1} = b_{\mu-1} = g_N^{(\mu-1)}-g_0^{(\mu-1)},\qquad\text{for}~\mu = 1,2,\ldots,\theta,
\end{equation*}
while vector $\displaystyle C_k\in\mathbb{C}^\theta_\text{c}$, which is entirely known for any given set of data, is defined by
\begin{equation*}
	(C_k)_{\mu,1} = 
    \begin{cases}
		-J_{0,k}f_{0,k}, & \text{if}\ \mu=1,\\ 
        0, & \text{if}\ 1<\mu\leq\theta,
	\end{cases}
\end{equation*}
where $f_{0,k}$ is the discrete Fourier transform of the sequence $g_j$ for $j=0,1,\ldots,N-1$. We therefore have
\begin{equation}
    f_{0,k}=\sum_{j=0}^{N-1}g_{j}\exp(-\mathrm{i}2\pi kj/N).
    \label{Eq_DFT}
\end{equation}

As $F_{\theta,k}$ is used to obtain the numerical derivatives, which are then used to construct the spline function, a good approximation of $F_{\theta,k}$ is needed. To achieve this, the accuracy of the boundary conditions vector $B$ plays a crucial role. A method to approximate this vector is presented in \cite{beaudoin_tutorial/article_1998,beaudoin_new_2003}. While it may yield interesting results in some cases, the method lacks robustness and can lead to inaccurate approximations of the boundary conditions when the parameters needed for its computation are not chosen adequately. Unfortunately, no general method to determine these parameters is currently suggested and the accuracy of the approximation established in \cite{beaudoin_tutorial/article_1998,beaudoin_new_2003} is strongly dependent on the information given at the nodes located near the center of the interval $[0,T]$. For these reasons, we present in this paper two algorithms for the computation of the boundary conditions. These two algorithms, which are presented in Section \ref{sec:AlgBC}, are independent of parameters and only the known values of a function at the interpolation nodes are needed as input.

Once the boundary conditions have been approximated, Equation \eqref{Eq_SystMatricielFb} can then be solved to determinate the discrete Fourier transform, $f_{\mu,k,\theta}$, for $\mu=1, 2, \ldots, \theta$. 

The approximation of the derivatives of $g$ at the interpolation nodes are then obtained by using the inverse discrete Fourier transform :
\begin{equation}
	g^{(\mu)}_{j,\theta}= \frac{1}{N}\sum_{k=0}^{N-1}f_{\mu,k,\theta}\exp(\mathrm{i}2\pi kj/N)
    \label{Eq_iDFT_theta}
\end{equation}
for $j=0,1,\ldots,N-1$ and $\mu=1,2,\ldots,\theta$.

\section{Approximation of the boundary conditions} 
\label{sec:AlgBC}

As previously stated, when the boundary conditions vector $B$ is properly computed, the interpolation method developed by Beaudoin and Beauchemin offers an accurate approximation of the numerical derivatives of a digitized function at the interpolation nodes. Moreover, using accurate approximations of the boundary conditions to solve Equation \eqref{Eq_SystMatricielFb} improves the accuracy of the overall interpolation as the degree of the spline $\theta$ is increased. We are therefore interested in determining two robust algorithms, independent of parameters such as $\theta$, that will accurately approximate the boundary conditions vector from any given set of $N+1$ equally spaced values. To do this, let us first note that, for any set $\{g_0,g_1,\ldots,g_N\}$ of $N+1$ equally spaced values of a digitized function $g$, $b_0$ is known since $b_0=g_N-g_0$. The values that need to be approximated are $b_1,\ldots,b_{\theta-1}$. Let us define the column vector $X\in\mathbb{R}_\text{c}^{\theta-1}$ by
\begin{equation*}
	X=\begin{bmatrix}b_1 \\ b_2 \\ \vdots \\ b_{\theta-1}\end{bmatrix}.
\end{equation*}
It then follows that
\begin{equation*}
	B=\begin{bmatrix}b_0 \\ X \end{bmatrix}.
\end{equation*}

To approximate this vector, our approach minimizes the square of the $L^2$-norm of the difference between two splines of consecutive degree (i.e. the difference between a spline of degree $\theta-1$ and degree $\theta$), which leads to the minimization of the following expression in terms of $X$
\begin{equation}
	\left\lVert g_\theta(t)-g_{\theta-1}(t)\right\rVert_{L^2}^2=\sum_{j=0}^{N-1}\left\lVert\left[g_\theta\right]_j(t)-\left[g_{\theta-1}\right]_j(t)\right\rVert_{L^2}^2.
    \label{Eq_ErrorNorm}
\end{equation}

\noindent Let us note that for any given value $\theta\in\mathbb{N}$, either $\theta-1$ or $\theta$ will be even, which means that $g_{\theta-1}$ or $g_\theta$ cannot be computed if $N$ is chosen as an even number. This implies that the minimization problem
\begin{equation}
	\underset{X\in\mathbb{R}_\text{c}^{\theta-1}}{\text{minimize}} \left\lVert g_\theta(t)-g_{\theta-1}(t)\right\rVert_{L^2}^2,
	\label{Eq_Minimization1}
\end{equation}
can only be solved for odd values of $N$. However, as we will see in the following section, \eqref{Eq_Minimization1} can be modified to obtain a minimization problem that works for all values of $N$ when $\theta$ is odd.

\noindent By expanding Equation \eqref{Eq_ErrorNorm}, we get
\begin{equation}
	\left\lVert g_\theta(t)-g_{\theta-1}(t)\right\rVert_{L^2}^2=\sum_{j=0}^{N-1}\int_{t_j}^{t_{j+1}}\left(\sum_{\alpha=1}^{\theta-1}\frac{(t-t_j)^\alpha}{\alpha!}\left(g_{j,\theta}^{(\alpha)}-g_{j,\theta-1}^{(\alpha)}\right)+\frac{(t-t_j)^\theta}{\theta!}g_{j,\theta}^{(\theta)}\right)^2\ud t,
    \label{Eq_normeDeLErreur}
\end{equation}
where $\displaystyle g_{j,\theta-1}^{(\alpha)}=\left[g^{(\alpha)}_{\theta-1}\right]_j(t_j)$ and $\displaystyle g_{j,\theta}^{(\alpha)}=\left[g^{(\alpha)}_\theta\right]_j(t_j)$.

From Equation \eqref{Eq_normeDeLErreur}, we propose two different methods for the approximation of the boundary conditions vector $B$. In the following, the dependence of $g_{\theta}^{(\beta)}$ on $B$ is omitted to simplify the notation.

\subsection{Method 1 - Minimizing the square of the $\theta$-th derivative of the spline of degree $\theta$}
\label{Section_BCfromApproxMin}
The method developed in this section is based on the assumption that the terms $g_{j,\theta}^{(\alpha)}-g_{j,\theta-1}^{(\alpha)}$ in Equation \eqref{Eq_normeDeLErreur}, for $\alpha = 1,2,\ldots,\theta-1$, can be neglected. By neglecting these terms, the equation no longer depends simultaneously on $\theta-1$ and $\theta$, which means that the equation can be evaluated for all values of $N$ when $\theta$ is odd. The singularity does however remain when $N$ and $\theta$ are both chosen as even integers. Under this assumption, the right-hand side of Equation \eqref{Eq_normeDeLErreur} becomes
\begin{equation*}
	\sum_{j=0}^{N-1}\int_{t_j}^{t_{j+1}}\left(\frac{(t-t_j)^\theta}{\theta!}g_{j,\theta}^{(\theta)}\right)^2\ud t,
\end{equation*}
which can be simplified to
\begin{equation}
	\sum_{j=0}^{N-1}\int_{t_j}^{t_{j+1}}\left(\frac{(t-t_j)^\theta}{\theta!}g_{j,\theta}^{(\theta)}\right)^2\ud t=\frac{(\Delta t)^{2\theta+1}}{(2\theta+1)(\theta!)^2}\sum_{j=0}^{N-1}\left(g_{j,\theta}^{(\theta)}\right)^2.
\label{Eq_Error_between_two_cons_orders}
\end{equation}
By minimizing \eqref{Eq_Error_between_two_cons_orders} in terms of $X$, the goal is to get an accurate approximation of the boundary conditions (i.e. $b_n$ for $n=1,\ldots,\theta-1$). This minimization problem can be written as
\begin{equation}
\begin{aligned}
& \underset{X\in\mathbb{R}_\text{c}^{\theta-1}}{\text{minimize}} && \left\lVert g_\theta(t)-g_{\theta-1}(t)\right\rVert_{L^2}^2 \approx\underset{X\in\mathbb{R}_\text{c}^{\theta-1}}{\text{minimize}}
& & \sum_{j=0}^{N-1}\left(g_{j,\theta}^{(\theta)}\right)^2,
\end{aligned}
\label{Eq_Minimization}
\end{equation}
or as
\begin{equation}
	\underset{X\in\mathbb{R}_\text{c}^{\theta-1}}{\text{minimize}}\int_0^T \left(g_\theta^{(\theta)}(t)\right)^2\ud t,
    \label{Eq_Minimization_integrale}
\end{equation}
since $g_\theta^{(\theta)}$ is constant on every subinterval $[t_j,t_{j+1})$.

\noindent The advantage of using \eqref{Eq_Minimization_integrale} versus \eqref{Eq_ErrorNorm} is that it can used for any integer $N$ when $\theta$ is odd.

Let us now develop the procedure for the first method.
By rewriting Equation \eqref{Eq_SystMatricielFb} as
\begin{equation*}
	F_{\theta,k} = M_{\theta,k}^{-1}(B+C_k),
\end{equation*}
the elements of vector $F_{\theta,k}$ are computed using the following equation
\begin{equation}
	f_{\alpha,k,\theta} = \text{row}_\alpha\left(M_{\theta,k}^{-1}\right)\cdot B - \left(M_{\theta,k}^{-1}\right)_{\alpha,1}J_{0,k}f_{0,k}, \qquad \text{for~}\alpha=1,2\ldots,\theta,
    \label{Eq_DFT_theta}
\end{equation}
where, from Equation \eqref{Eq_MatrixInverseMtheta}, the scalar $\left(M_{\theta,k}^{-1}\right)_{\alpha,1}\in\mathbb{C}$ is computed as follows
\begin{equation*}
    \left(M_{\theta,k}^{-1}\right)_{\alpha,1}=(-J_{0,k})^{\alpha-1}\frac{\det(M_{\theta-\alpha,k})}{\det(M_{\theta,k})},\qquad \text{for~}\alpha=1,2\ldots,\theta.
\end{equation*}
By replacing Equation \eqref{Eq_DFT_theta} in Equation \eqref{Eq_iDFT_theta} for $\alpha=\theta$, it follows that
\begin{equation*}
	g^{(\theta)}_{j,\theta} = \Delta_{\theta,j}\cdot B - \sigma_{\theta,j},
\end{equation*}
where $\displaystyle \Delta_{\theta,j}\in\mathbb{C}^\theta_\text{r}$ is defined by
\begin{equation}
	(\Delta_{\theta,j})_{1,\nu} = \delta_{\nu-1,j,\theta} = \frac{1}{N}\sum_{k=0}^{N-1}(-J_{0,k})^{\theta-\nu}\frac{\det(M_{\nu-1,k})}{\det(M_{\theta,k})}\exp(\mathrm{i}2\pi kj/N),\qquad\text{for }\nu=1,2,\ldots,\theta, 
         \label{Eq_DeltaTheta}
\end{equation}
and $\displaystyle\sigma_{\theta,j}\in\mathbb{C}$ is defined by 
\begin{equation}
	\sigma_{\theta,j} = \frac{(-1)^{\theta-1}}{N}\sum_{k=0}^{N-1}\frac{(J_{0,k})^\theta}{\det(M_{\theta,k})} f_{0,k}\exp(\mathrm{i}2\pi kj/N),
	\label{Eq_SigmaDelta}
\end{equation}
for all $j=0,1,\ldots,N-1$.

\noindent Let $\displaystyle h:\mathbb{R}_\text{c}^\theta\rightarrow\mathbb{R}$ be defined by $\displaystyle h(B) = \sum_{j=1}^{N-1}\left(g_{j,\theta}^{(\theta)}\right)^2$. It follows that
\begin{equation*}
	h(B) = \sum_{j=0}^{N-1}(\Delta_{\theta,j}\cdot B-\sigma_{\theta,j})^2.
    \end{equation*}
As the minimization of $h$ in terms of $X$ is a convex problem, the minimum can be obtained by solving
\begin{equation}
	\nabla_{X} h(B) = \pmb{0}_{(\theta-1)\times1},
    \label{Eq_gradient_de_A}
\end{equation}
where $\pmb{0}_{(\theta-1)\times1}$ is the null column vector of dimension $\theta-1$.

\noindent By expanding \eqref{Eq_gradient_de_A}, we get 
\begin{equation*}
	\sum_{\alpha=0}^{\theta-1}b_\alpha\sum_{j=0}^{N-1} \delta_{n,j,\theta}\delta_{\alpha,j,\theta}=\sum_{j=0}^{N-1}\sigma_{\theta,j}\delta_{n,j,\theta},
\end{equation*}
for $n=1,\dots,\theta-1$, which can be rewritten as
\begin{equation}
	\sum_{\alpha=1}^{\theta-1}b_\alpha\sum_{j=0}^{N-1}\delta_{n,j,\theta}\delta_{\alpha,j,\theta}=\sum_{j=0}^{N-1}(\sigma_{\theta,j}\delta_{n,j,\theta}-b_0\delta_{n,j,\theta}\delta_{0,j,\theta}),
    \label{Eq_ApproxNob0BC}
\end{equation}
since $b_0$ is known. 

Equation \eqref{Eq_ApproxNob0BC} is a system of $\theta-1$ linear equations, which can be written as
\begin{equation}
	\Gamma X = \Sigma,
    \label{Eq_SysMatriciel_CF_sans_b0}
\end{equation}
where $\Gamma$ is a square matrix of dimension $\theta-1$ defined by
\begin{equation}
	(\Gamma)_{\mu,\nu} = \sum_{j=0}^{N-1} \delta_{\mu,j,\theta}\delta_{\nu,j,\theta}, 
	\label{Eq_MatrixGamma}
\end{equation}
and $\displaystyle \Sigma\in\mathbb{C}_\text{c}^{\theta-1}$ is defined by
\begin{equation}
	(\Sigma)_{\mu,1} = \sum_{j=0}^{N-1} \left(\sigma_{\theta,j}\delta_{\mu,j,\theta}-b_0\delta_{\mu,j,\theta}\delta_{0,j,\theta}\right),
	\label{Eq_MatrixSigma}
\end{equation}
for $\mu,\nu=1,2,\ldots,\theta-1$.

\noindent By solving system \eqref{Eq_SysMatriciel_CF_sans_b0}, we obtain the desired approximation of $X$. To solve the system, different methods can be used, but in this paper, we use the LU decomposition method to avoid a direct computation of matrix $\Gamma^{-1}$. 
Numerical results are shown in Section \ref{Sec:Numerical_results} and illustrate the good performance of our approach.

The numerical implementation of Method 1 to approximate the boundary conditions can be summarized as follows :
\begin{itemize}
    \item[$\bullet$] Computation of $f_{0,k}$ for $k=0,1,\ldots,N-1$ using Equation \eqref{Eq_DFT};
    \item[$\bullet$] Computation of $\Delta_{\theta,j}$ for $j=0,1,\ldots,N-1$ using Equation \eqref{Eq_DeltaTheta};
    \item[$\bullet$] Computation of $\sigma_{\theta,j}$ for $j=0,1,\ldots,N-1$ using Equation \eqref{Eq_SigmaDelta};
    \item[$\bullet$] Computation of $\Gamma$ using Equation \eqref{Eq_MatrixGamma};
    \item[$\bullet$] Computation of $\Sigma$ using Equation \eqref{Eq_MatrixSigma}; and
    \item[$\bullet$] Computation of $X$ by solving System \eqref{Eq_SysMatriciel_CF_sans_b0}.
\end{itemize}

\subsection{Method 2 - Minimizing the difference between splines of consecutive degrees}
\label{Section_BCfromEM}

In the previous method, we assumed that the terms $g_{j,\theta}^{(\alpha)}-g_{j,\theta-1}^{(\alpha)}$ could be neglected for \newline $\alpha=1,2,\ldots,\theta-1$. Without this assumption, Equation \eqref{Eq_normeDeLErreur} can be written as
\begin{align*}
\begin{split}
    \left\lVert g_\theta(t)-g_{\theta-1}(t)\right\rVert_{L^2}^2&=\sum_{j=0}^{N-1}\int_{t_j}^{t_{j+1}}\left(\sum_{\alpha=1}^\theta\frac{(t-t_j)^\alpha}{\alpha!}G_{j,\theta}^{(\alpha)}\right)^2\ud t\\
    &=\sum_{j=0}^{N-1}\sum_{\alpha=1}^\theta\int_{t_j}^{t_{j+1}}\left\{\left(\frac{(t-t_j)^\alpha}{\alpha!}G_{j,\theta}^{(\alpha)}\right)^2+2\sum_{\beta=1}^{\alpha-1}\frac{(t-t_j)^{\alpha+\beta}}{\alpha!\beta!}G_{j,\theta}^{(\alpha)}G_{j,\theta}^{(\beta)}\right\}\ud t\\
    &=\sum_{j=0}^{N-1}\sum_{\alpha=1}^\theta\left\{\zeta(\alpha,\alpha)\left(G_{j,\theta}^{(\alpha)}\right)^2+2\sum_{\beta=1}^{\alpha-1}\zeta(\alpha,\beta)G_{j,\theta}^{(\alpha)}G_{j,\theta}^{(\beta)}\right\},
\end{split}
\end{align*}
where
\begin{equation}
	\zeta(\alpha,\beta)=\frac{(\Delta t)^{\alpha+\beta+1}}{(\alpha+\beta+1)\alpha!\beta!},
	\label{Eq_zeta}
\end{equation}
and
  \begin{equation}
    G_{j,\theta}^{(\alpha)} =
    \begin{cases*}
      g_{j,\theta}^{(\alpha)}-g_{j,\theta-1}^{(\alpha)}, & if $1\leq\alpha\leq\theta-1$, \\
      g_{j,\theta}^{(\theta)},        & if $\alpha=\theta$.
    \end{cases*}
    \label{Eq_GjAlpha_genCase}
  \end{equation}

\noindent If $N$ is odd, both $M_{\theta-1,k}$ and $M_{\theta,k}$ are invertible \cite{pepin_new_2019}. Let us define the row vector $\displaystyle D_{\alpha,k}\in\mathbb{C}_\text{r}^{\theta}$ by
\begin{equation*}
    (D_{\alpha,k})_{1,\nu} = d_{\nu-1,k,\alpha} =
    \begin{cases*}
      \left(M_{\theta-1,k}^{-1}\right)_{\alpha,\nu}, & if $1\leq\nu\leq\theta-1$, \\
      0,        & if $\nu=\theta$,
    \end{cases*}
  \end{equation*}
for $\alpha = 1,2,\ldots,\theta-1$. As $\displaystyle M_{\theta-1,k}^{-1}$ is of dimension $\theta-1$,  $D_{\theta,k}\in\mathbb{C}_\text{r}^{\theta}$ will be defined by convention as $D_{\theta,k}=\pmb{0}_{1\times\theta}$ .

\noindent Similarly, the row vectors $\displaystyle \Omega_{\alpha,j}\in\mathbb{C}_\text{r}^{\theta}$ and $\displaystyle \Delta_{\alpha,j}\in\mathbb{C}_\text{r}^{\theta}$ will be defined by

\begin{equation}
	\Omega_{\alpha,j} = \frac{1}{N}\sum_{k=0}^{N-1}D_{\alpha,k}\exp(\mathrm{i}2\pi kj/N) = \begin{bmatrix}
           \omega_{0,j,\alpha} & \omega_{1,j,\alpha} & \cdots & \omega_{\theta-2,j,\alpha} & 0
         \end{bmatrix},
         \label{Eq_OmegaGenCase}
\end{equation}
and
\begin{equation}
	\Delta_{\alpha,j} = \frac{1}{N}\sum_{k=0}^{N-1}\text{row}_\alpha\left(M_{\theta,k}^{-1}\right)\exp(\mathrm{i}2\pi kj/N) = \begin{bmatrix}
           \delta_{0,j,\alpha} & \delta_{1,j,\alpha} & \cdots & \delta_{\theta-1,j,\alpha}
         \end{bmatrix},
         \label{Eq_DeltaGenCase}
\end{equation}
for $\alpha=1,2,\ldots,\theta$. When $\alpha=\theta$, it follows that $\Omega_{\theta,j}=\pmb{0}_{1\times\theta}$. Explicitly, Equations \eqref{Eq_OmegaGenCase} and \eqref{Eq_DeltaGenCase} can be written as
\begin{equation}
	(\Omega_{\alpha,j})_{1,\nu} = \omega_{\nu-1,j,\alpha} = \begin{cases*}
      \displaystyle\frac{1}{N}\sum_{k=0}^{N-1}d_{\nu-1,k,\alpha}\exp(\mathrm{i}2\pi kj/N), & if $1\leq\alpha\leq\theta-1$, \\
      0,        & if $\alpha=\theta$,
    \end{cases*}
    \label{Eq_Gammaalpha}
\end{equation}
and 
\begin{equation}
	(\Delta_{\alpha,j})_{1,\nu} = \delta_{\nu-1,j,\alpha} = \frac{1}{N}\sum_{k=0}^{N-1}\left(M_{\theta,k}^{-1}\right)_{\alpha,\nu}\exp(\mathrm{i}2\pi kj/N),
	\label{Eq_Deltaalpha}
\end{equation}
for $\alpha,\nu = 1,2,\ldots,\theta$. Let us note that equations \eqref{Eq_Gammaalpha} and \eqref{Eq_Deltaalpha} are both inverse discrete Fourier transforms which can be computed using a FFT algorithm. 

\noindent Let us also define $\eta_{\alpha,j}\in\mathbb{C}$ and $\sigma_{\alpha,j}\in\mathbb{C}$ by
\begin{equation}
	  \eta_{\alpha,j} = \begin{cases*}
      \displaystyle\frac{(-1)^{\alpha-1}}{N}\sum_{k=0}^{N-1}(J_{0,k})^\alpha\frac{\det(M_{\theta-1-\alpha,k})}{\det(M_{\theta-1,k})}f_{0,k}\exp(\mathrm{i}2\pi kj/N), & if $1\leq\alpha\leq\theta-1$, \\
      0,        & if $\alpha=\theta$,
    \end{cases*}
    \label{Eq_Etaalpha}
\end{equation}
and
\begin{equation}
	\sigma_{\alpha,j} = \frac{(-1)^{\alpha-1}}{N}\sum_{k=0}^{N-1}(J_{0,k})^\alpha\frac{\det(M_{\theta-\alpha,k})}{\det(M_{\theta,k})}f_{0,k}\exp(\mathrm{i}2\pi kj/N),
	\label{Eq_Sigmaalpha}
\end{equation}
for $\alpha=1,2,\ldots,\theta$. 
\noindent Equations \eqref{Eq_Etaalpha} and \eqref{Eq_Sigmaalpha} are also inverse discrete Fourier transforms. A FFT algorithm can therefore be used to compute them.

\noindent The numerical derivatives of $g$ of an approximation order $\theta-1$ can then be written as
\begin{equation}
	g_{j,\theta-1}^{(\alpha)} = \Omega_{\alpha,j}\cdot B-\eta_{\alpha,j},
    \label{Eq_GenCase_gj_thetam1}
\end{equation}
while those of order $\theta$ can be written as
\begin{equation}
	g_{j,\theta}^{(\alpha)} = \Delta_{\alpha,j}\cdot B-\sigma_{\alpha,j}.
     \label{Eq_GenCase_gj_theta}
\end{equation}

\noindent Furthermore, by defining the row vector $\displaystyle \Psi_{\alpha,j}\in\mathbb{C}_\text{r}^\theta$ by
\begin{equation}
    \Psi_{\alpha,j} = \Delta_{\alpha,j}-\Omega_{\alpha,j} = \begin{bmatrix}
           \psi_{0,j,\alpha} & \psi_{1,j,\alpha} & \cdots & \psi_{\theta-1,j,\alpha}
         \end{bmatrix},
         \label{Eq_Psialpha}
\end{equation}
where
\begin{equation*}
    (\Psi_{\alpha,j})_{1,\nu} = \psi_{\nu-1,j,\alpha} = \delta_{\nu-1,j,\alpha}-\omega_{\nu-1,j,\alpha},\qquad\text{for}~\alpha,\nu=1,2,\ldots,\theta,
\end{equation*}
and $\phi_{\alpha,j}\in\mathbb{C}$ by
\begin{equation}
    \phi_{\alpha,j}=\sigma_{\alpha,j}-\eta_{\alpha,j},\qquad\text{for~}\alpha=1,2,\ldots,\theta,
    \label{Eq_Phialpha}
\end{equation}
Equation \eqref{Eq_GjAlpha_genCase} can be rewritten as
\begin{align*}
\begin{split}
	G^{(\alpha)}_{j,\theta}&=\left(\Delta_{\alpha,j}-\Omega_{\alpha,j}\right)\cdot B-(\sigma_{\alpha,j}-\eta_{\alpha,j})\\
  &=\Psi_{\alpha,j}\cdot B-\phi_{\alpha,j}.
\end{split}
\end{align*}
by using Equations \eqref{Eq_GenCase_gj_thetam1} and \eqref{Eq_GenCase_gj_theta}.

\noindent Since the function $h^*: \mathbb{R}_\text{c} ^\theta\rightarrow\mathbb{R}$ defined by $h^*(B)=\left\lVert g_\theta(t)-g_{\theta-1}(t)\right\rVert_{L^2}^2$ is a convex function in terms of $X$, its minimum value is obtained by solving
\begin{equation}
	\nabla_{X} h^*(B)=\pmb{0}_{(\theta-1)\times1},
    \label{Eq_Gradient_of_h_star}
\end{equation} 
where $\nabla_{X}$ is the gradient with respect to $X$. It can then be shown that Equation \eqref{Eq_Gradient_of_h_star} is equivalent to
\begin{equation}
	\sum_{j=0}^{N-1}\sum_{\alpha=1}^\theta\left\{\zeta(\alpha,\alpha)\frac{\partial}{\partial b_n}\left(G_{j,\theta}^{(\alpha)}\right)^2+2\sum_{\beta=1}^{\alpha-1}\zeta(\alpha,\beta)\frac{\partial}{\partial b_n}\left(G_{j,\theta}^{(\alpha)}G_{j,\theta}^{(\beta)}\right)\right\}=0,
    \label{Eq_GenCase_DerOfNorm}
\end{equation}
for $n=1,\ldots,\theta-1$.

\noindent One can easily verify that
\begin{equation}
	\frac{\partial}{\partial b_n}\left(G_{j,\theta}^{(\alpha)}\right)^2 = 2\psi_{n,j,\alpha}\left(\sum_{a=0}^{\theta-1}b_a\psi_{a,j,\alpha}-\phi_{\alpha,j}\right),
    \label{Eq_GenCase_PartialDerGjsquared}
\end{equation}
and
\begin{equation}
	\frac{\partial}{\partial b_n}\left(G_{j,\theta}^{(\alpha)}G_{j,\theta}^{(\beta)}\right) = \sum_{a=0}^{\theta-1}b_a\left(\psi_{a,j,\alpha}\psi_{n,j,\beta}+\psi_{n,j,\alpha}\psi_{a,j,\beta}\right)-\psi_{n,j,\alpha}\phi_{\beta,j}-\psi_{n,j,\beta}\phi_{\alpha,j},
    \label{Eq_GenCase_PartialDerGjGj}
\end{equation}
for $n=1,\ldots,\theta-1$.

\noindent By replacing Equations \eqref{Eq_GenCase_PartialDerGjsquared} and \eqref{Eq_GenCase_PartialDerGjGj} in \eqref{Eq_GenCase_DerOfNorm}, it then follows that
\begin{align*}
\begin{split}
	\sum_{a=0}^{\theta-1}b_a\sum_{j=0}^{N-1}\sum_{\alpha=1}^{\theta}&\left\{\zeta(\alpha,\alpha)\psi_{n,j,\alpha}\psi_{a,j,\alpha}+\sum_{\beta=1}^{\alpha-1}\zeta(\alpha,\beta)(\psi_{a,j,\alpha}\psi_{n,j,\beta}+\psi_{n,j,\alpha}\psi_{a,j,\beta})\right\}\\
    &=\sum_{j=0}^{N-1}\sum_{\alpha=1}^{\theta}\left\{\zeta(\alpha,\alpha)\psi_{n,j,\alpha}\phi_{\alpha,j}+\sum_{\beta=1}^{\alpha-1}\zeta(\alpha,\beta)(\psi_{n,j,\alpha}\phi_{\beta,j}+\psi_{n,j,\beta}\phi_{\alpha,j})\right\},
\end{split}
\end{align*}
for $n=1,\ldots,\theta-1$.

\noindent Since $b_0$ is known, we can write the last equation as
\begin{align*}
\begin{split}
	\sum_{a=1}^{\theta-1}b_a\sum_{j=0}^{N-1}\sum_{\alpha=1}^{\theta}&\left\{\zeta(\alpha,\alpha)\psi_{n,j,\alpha}\psi_{a,j,\alpha}+\sum_{\beta=1}^{\alpha-1}\zeta(\alpha,\beta)(\psi_{a,j,\alpha}\psi_{n,j,\beta}+\psi_{n,j,\alpha}\psi_{a,j,\beta})\right\}\\
    &=\sum_{j=0}^{N-1}\sum_{\alpha=1}^{\theta}\left\{\zeta(\alpha,\alpha)\psi_{n,j,\alpha}\phi_{\alpha,j}+\sum_{\beta=1}^{\alpha-1}\zeta(\alpha,\beta)(\psi_{n,j,\alpha}\phi_{\beta,j}+\psi_{n,j,\beta}\phi_{\alpha,j})-\rho_{n,j,\alpha}\right\},
\end{split}
\end{align*}
for $n=1,2,\ldots,\theta-1$ where
\begin{equation}
	\rho_{n,j,\alpha} = b_0\left\{\zeta(\alpha,\alpha)\psi_{n,j,\alpha}\psi_{0,j,\alpha}+\sum_{\beta=1}^{\alpha-1}\zeta(\alpha,\beta)(\psi_{0,j,\alpha}\psi_{n,j,\beta}+\psi_{n,j,\alpha}\psi_{0,j,\beta})\right\}.
	\label{Eq_Rhoalpha}
\end{equation}
This system of linear equations can be written in matrix form as
\begin{equation}
	\Lambda X=\Pi,
	\label{Eq_M2SystemToSolve}
\end{equation}
where $\Lambda$ is a square matrix of dimension $\theta-1$ with general term
\begin{equation}
(\Lambda)_{\mu,\nu} = \sum_{j=0}^{N-1}\sum_{\alpha=1}^{\theta}\left\{\zeta(\alpha,\alpha)\psi_{\nu,j,\alpha}\psi_{\mu,j,\alpha}+\sum_{\beta=1}^{\alpha-1}\zeta(\alpha,\beta)(\psi_{\mu,j,\alpha}\psi_{\nu,j,\beta}+\psi_{\nu,j,\alpha}\psi_{\mu,j,\beta})\right\},
\label{Eq_Lambda}
\end{equation}
and $\Pi\in\mathbb{C}^{\theta-1}_\text{c}$ is defined by
\begin{equation}
	(\Pi)_{\mu,1} = \sum_{j=0}^{N-1}\sum_{\alpha=1}^{\theta}\left\{\zeta(\alpha,\alpha)\psi_{\mu,j,\alpha}\phi_{\alpha,j}+\sum_{\beta=1}^{\alpha-1}\zeta(\alpha,\beta)(\psi_{\mu,j,\alpha}\phi_{\beta,j}+\psi_{\mu,j,\beta}\phi_{\alpha,j})-\rho_{\mu,j,\alpha}\right\},
    \label{Eq_Pi}
\end{equation}
for $\mu,\nu=1,2,\ldots,\theta-1$. 

Method 2 therefore also requires solving a linear system of equations to approximate the boundary conditions vector $X$. However, in this case, System \eqref{Eq_M2SystemToSolve} can only be solved for odd integer values of $N$ due to the nature of the minimization problem.
The numerical implementation of this method can be summarized as follows :

\begin{itemize}
    \item[$\bullet$] Computation of $f_{0,k}$ for $k=0,1,\ldots,N-1$ using Equation \eqref{Eq_DFT};
    \item[$\bullet$] Computation of $\Psi_{\alpha,j}$ for $\alpha=1,2,\ldots,\theta$ and $j=0,1,\ldots,N-1$ using Equation \eqref{Eq_Psialpha};
    \item[$\bullet$] Computation of $\phi_{\alpha,j}$ for $\alpha=1,2,\ldots,\theta$ and $j=0,1,\ldots,N-1$ using Equation \eqref{Eq_Phialpha};
    \item[$\bullet$] Computation of $\rho_{n,j,\alpha}$ for $n=1,2,\ldots,\theta-1$, $\alpha=1,2,\ldots,\theta$ and $j=0,1,\ldots,N-1$ using Equation \eqref{Eq_Rhoalpha};
    \item[$\bullet$] Computation of $\Lambda$ using Equation \eqref{Eq_Lambda};
    \item[$\bullet$] Computation of $\Pi$ using Equation \eqref{Eq_Pi}; and
    \item[$\bullet$] Computation of $X$ by solving System \eqref{Eq_M2SystemToSolve}.
\end{itemize}

\section{Numerical results}
\label{Sec:Numerical_results}
In order to verify the effectiveness of the two proposed methods, they will be compared with the traditional cubic spline interpolation method, for both natural and not-a-knot boundary conditions (see \cite{de_boor_practical_1978} for more information on these boundary conditions). The methods will then be used to compute spline functions of degrees $4$, $5$ and $11$ in order to study the gain in accuracy obtained by increasing $\theta$. The accuracy of the methods will be determined by computing an average and maximum error between the analytic functions and the spline functions. As for the comparisons, four analytic functions, defined on different intervals, will be considered. As mentioned in Section \ref{Section_BCfromEM}, even numbers of interpolation points (i.e. odd values of $N$) will be used to ensure that both $M_{\theta-1,k}$ and $M_{\theta,k}$ are invertible. 

\begin{figure}[htbp!] 
  \begin{subfigure}[t]{0.5\linewidth}
    \centering
    \includegraphics[width=\linewidth]{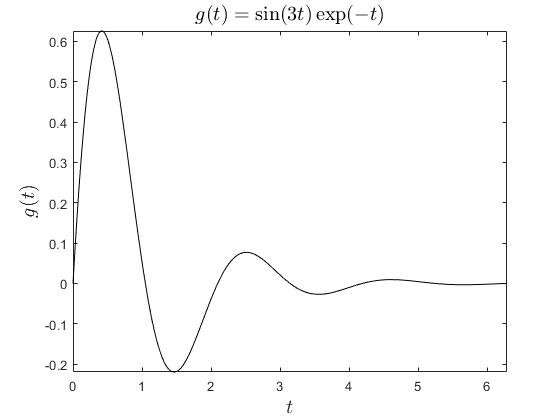} 
    \caption{Interval: $[0,2\pi]$.}
    \label{fig7:a} 
    \vspace{4ex}
  \end{subfigure}
  \begin{subfigure}[t]{0.5\linewidth}
    \centering
    \includegraphics[width=\linewidth]{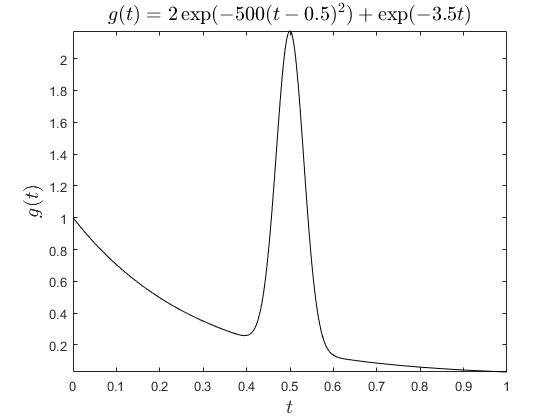} 
        \caption{Interval: $[0,1]$.} 
    \label{fig7:b} 
    \vspace{4ex}
  \end{subfigure} 
  \begin{subfigure}[t]{0.5\linewidth}
    \centering
    \includegraphics[width=\linewidth]{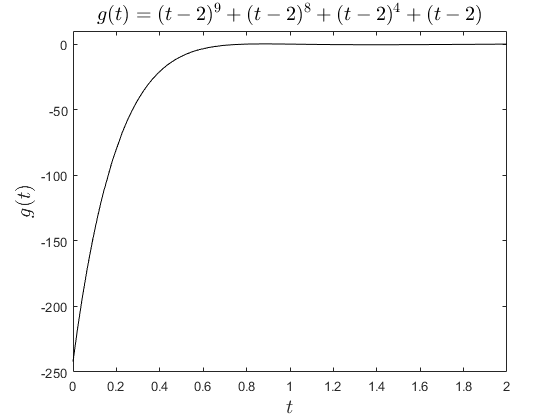} 
    \caption{Interval: $[0,2]$.}
    \label{fig7:c} 
  \end{subfigure}
  \begin{subfigure}[t]{0.5\linewidth}
    \centering
    \includegraphics[width=\linewidth]{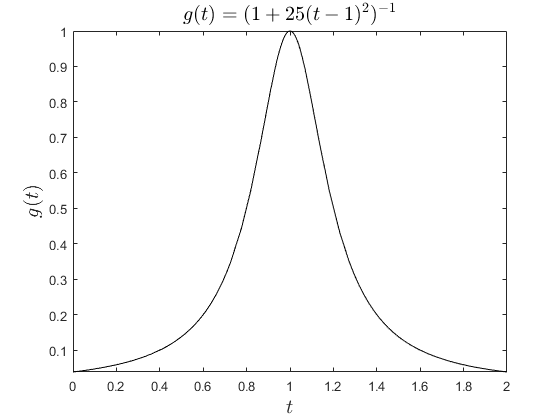} 
    \caption{Interval: $[0,2]$.}
    \label{fig7:d} 
  \end{subfigure} 
  \caption{Analytic functions to interpolate on $[0,T]$.}
  \label{fig7} 
\end{figure}

\noindent Let us now consider the following four analytic functions:
\begin{enumerate}
\item[$\bullet$] $g(t) = \sin(3t)\exp(-t)$ on the interval $[0,2\pi]$;
\item[$\bullet$] $g(t)=2\exp(-500(t-0.5)^2)+\exp(-3.5t)$ on the interval $[0,1]$;
\item[$\bullet$] $g(t)=(t-2)^9+(t-2)^8+(t-2)^4+(t-2)$ on the interval $[0,2]$;
\item[$\bullet$] $g(t)=\displaystyle(1+25(t-1)^2)^{-1}$ on the interval $[0,2]$.
\end{enumerate}
These functions are illustrated in Figure \ref{fig7}. They may seem easy to interpolate, but in reality, some are quite difficult, especially when using low values of $N$, due to their high peaks or oscillations. 

The spline functions will be computed using the method described in Section \ref{sec:MathBack} and the boundary conditions will be approximated using the methods presented in Section \ref{sec:AlgBC}.   
Let us note that in all cases, $b_0$ is computed analytically from the initial data.  Let us also note that other methods, such as setting all the boundary conditions to $0$ or finding the boundary conditions vector $X$ that minimizes $\left\lVert g_\theta(t)-g_{\theta-2}(t)\right\rVert_{L^2}^2$, were also tested. As these methods were less efficient than the proposed methods, they will not be presented.

\subsection{Error estimation}
The spline approximation $g_\theta$ of $g$ on the interval $[0,T]$ is defined by Equation \eqref{Polynome_fct_g}. To approximate the error between $g_\theta$ and $g$ on the interval $[0,T]$, both functions will be evaluated at $\lambda>1$ points on the subinterval $[t_j,t_{j+1}]$, $j=0,1,\ldots,N-1$, which leads to a total computation of  $\lambda N$ values. As for the maximum and average errors between the spline and the analytic function, they will be calculated, respectively, using the following formulas
\begin{equation*}
	E_\theta^{\text{max}} = \max_{\substack{j=0,1,\ldots,N-1\\m=1,2,\ldots,\lambda-1\\}}|g(t_j+m\Delta t_s)-g_\theta(t_j+m\Delta t_s)|,
\end{equation*}
and
\begin{equation*}
	E_\theta^{\text{avg}}=\frac{1}{(\lambda-1)N}\sum_{\gamma=0}^{\lambda N-1}|g(\gamma\Delta t_s)-g_\theta(\gamma\Delta t_s)|,
\end{equation*}
with $$\displaystyle g_\theta(t_j+m\Delta t_s)=\sum_{p=0}^\theta\frac{(m\Delta t_s)^p}{p!}g^{(p)}_{j,\theta},$$ and $$\displaystyle\Delta t_s=\frac{\Delta t}{\lambda},$$ for $j=0,1,\ldots,N-1$ and $m=1,2,\ldots,\lambda-1$.

Since the value of $|g(t)-g_\theta(t)|$ is zero at the interpolation nodes $t=t_j$, $j=0,1,\ldots,N-1$, we divide by $(\lambda-1)N$ instead of $\lambda N$ in the calculation of the average error in order to get a better approximation.

\subsection{Error analysis of cubic splines}

In practice, cubic splines are the most popular choice  for piecewise interpolation of a smooth function due to their accuracy, low computational cost and simplicity. Traditionally, cubic splines have two free parameters, one at each boundary of the interval $[0,T]$. These parameters are called the boundary conditions. Generally, they correspond to the values of the first derivative of the digitized function at $t=0$ and $t=T$. Since these two free parameters are typically not known in practice, they must be approximated or set to a default value. To compare this standard method with our proposed methods, we will use the natural and not-a-knot boundary conditions. The results will be obtained using the \textit{csape}\footnote{The \textit{Curve Fitting Toolbox\texttrademark} is required to use the \textit{csape} function.} function from MATLAB\textsuperscript{\textregistered}.

Tables \ref{Table_theta3_fct1}, \ref{Table_theta3_fct2}, \ref{Table_theta3_fct3} and \ref{Table_theta3_fct4} present the results obtained for the four analytic functions considered in the case of cubic splines. In these tables, NS stands for the traditional cubic splines with natural boundary conditions while NAK stands for the traditional cubic splines with not-a-knot boundary conditions. As for Methods 1 and 2, they represent the interpolation method developed by Beaudoin and Beauchemin, combined with the two methods presented in Section \ref{sec:AlgBC} to compute the boundary conditions. Let us also note that the results were obtained using $\lambda=10$ and that the values in bold font show the methods that offer the best results.

From these results, it can be seen that the not-a-knot boundary conditions seems to perform best in terms of accuracy. When looking at the results obtained with our proposed methods for $\theta=3$, it can also be seen that they are very similar to those obtained with the not-a-knot boundary conditions. Our proposed methods therefore lead to very accurate results and perform better than the traditional cubic splines with natural boundary conditions.

\begin{table}[!htbp]
\centering
\begin{tabular}{|Sc|Sc|Sc|Sc|Sc|Sc|}
\cline{1-6}
\multicolumn{2}{|c|}{$\theta=3$} & \multicolumn{4}{c|}{Method}     \\ \hline
$N$                     & Error                                                    & Method 1                                                      & Method 2                                                                                                                                                                                & NS                                                                                                      & NAK                                                      \\ \hline
                        & $E_\theta^{\text{max}}$ & \begin{tabular}[c]{@{}c@{}}$2.43\times10^{-03}$\end{tabular} & \begin{tabular}[c]{@{}c@{}}\pmb{$1.44\times10^{-03}$}\end{tabular} &  { \begin{tabular}[c]{@{}c@{}}$1.31\times10^{-02}$\end{tabular}} & \begin{tabular}[c]{@{}l@{}}$3.59\times10^{-03}$\end{tabular} \\ \cline{2-6} 
\multirow{-2}{*}{$31$}  & $E_\theta^{\text{avg}}$  & \begin{tabular}[c]{@{}c@{}}$9.44\times10^{-05}$\end{tabular} & \begin{tabular}[c]{@{}c@{}}\pmb{$7.51\times10^{-05}$}\end{tabular} & { \begin{tabular}[c]{@{}c@{}}$4.03\times10^{-04}$\end{tabular}} & \begin{tabular}[c]{@{}l@{}}$1.28\times10^{-04}$\end{tabular} \\ \hline
                        & $E_\theta^{\text{max}}$ & \begin{tabular}[c]{@{}c@{}}$1.07\times10^{-04}$\end{tabular} & \begin{tabular}[c]{@{}c@{}}$\pmb{5.55\times10^{-06}}$\end{tabular} & { \begin{tabular}[c]{@{}c@{}}$1.15\times10^{-03}$\end{tabular}} & \begin{tabular}[c]{@{}l@{}}$3.92\times10^{-05}$\end{tabular} \\ \cline{2-6} 
\multirow{-2}{*}{$101$} & $E_\theta^{\text{avg}}$ & \begin{tabular}[c]{@{}c@{}}$1.19\times10^{-06}$\end{tabular} & \begin{tabular}[c]{@{}c@{}}$7.43\times10^{-07}$\end{tabular} & { \begin{tabular}[c]{@{}c@{}}$1.08\times10^{-05}$\end{tabular}} & \begin{tabular}[c]{@{}l@{}}\pmb{$5.65\times10^{-07}$}\end{tabular} \\ \hline
                        & $E_\theta^{\text{max}}$ & \begin{tabular}[c]{@{}c@{}}$9.79\times10^{-07}$\end{tabular} & \begin{tabular}[c]{@{}c@{}}$2.77\times10^{-07}$\end{tabular} & { \begin{tabular}[c]{@{}c@{}}$4.63\times10^{-05}$\end{tabular}} & \begin{tabular}[c]{@{}l@{}}\pmb{$6.65\times10^{-08}$}\end{tabular} \\ \cline{2-6} 
\multirow{-2}{*}{$501$} & $E_\theta^{\text{avg}}$ & \begin{tabular}[c]{@{}c@{}}$2.21\times10^{-09}$\end{tabular} & \begin{tabular}[c]{@{}c@{}}$1.32\times10^{-09}$\end{tabular} & { \begin{tabular}[c]{@{}c@{}}$8.69\times10^{-08}$\end{tabular}} & \begin{tabular}[c]{@{}l@{}}\pmb{$5.18\times10^{-10}$}\end{tabular} \\ \hline
\end{tabular}
\caption{Comparisons in the case of $\theta=3$ for the function $g(t)=\sin(3t)e^{-t}$ on the interval $[0,2\pi]$. }
\label{Table_theta3_fct1}
\end{table}

\begin{table}[!htbp]
\centering
\begin{tabular}{|Sc|Sc|Sc|Sc|Sc|Sc|}
\cline{1-6}
\multicolumn{2}{|c|}{$\theta=3$} & \multicolumn{4}{c|}{Method}     \\ \hline
$N$                     & Error                   & Method 1                                                      & Method 2                                                                                                                                                                             & NS                                                                                                      & NAK                                                      \\ \hline
                        & $E_\theta^{\text{max}}$ & \begin{tabular}[c]{@{}c@{}}$4.39\times10^{-02}$\end{tabular} & \begin{tabular}[c]{@{}c@{}}$4.39\times10^{-02}$\end{tabular} & { \begin{tabular}[c]{@{}c@{}}$4.39\times10^{-02}$\end{tabular}} & \begin{tabular}[c]{@{}l@{}}$\pmb{4.29\times10^{-02}}$\end{tabular} \\ \cline{2-6} 
\multirow{-2}{*}{$31$}  & $E_\theta^{\text{avg}}$ & \begin{tabular}[c]{@{}c@{}}$\pmb{2.40\times10^{-03}}$\end{tabular} & \begin{tabular}[c]{@{}c@{}}$2.77\times10^{-03}$\end{tabular} & { \begin{tabular}[c]{@{}c@{}}$2.42\times10^{-03}$\end{tabular}} & \begin{tabular}[c]{@{}l@{}}$\pmb{2.40\times10^{-03}}$\end{tabular} \\ \hline
                        & $E_\theta^{\text{max}}$ & \begin{tabular}[c]{@{}c@{}}$\pmb{1.70\times10^{-04}}$\end{tabular} & \begin{tabular}[c]{@{}c@{}}$\pmb{1.70\times10^{-04}}$\end{tabular} & \begin{tabular}[c]{@{}c@{}}$\pmb{1.70\times10^{-04}}$\end{tabular} & \begin{tabular}[c]{@{}l@{}}$\pmb{1.70\times10^{-04}}$\end{tabular} \\ \cline{2-6} 
\multirow{-2}{*}{$101$} & $E_\theta^{\text{avg}}$ & \begin{tabular}[c]{@{}c@{}}$7.59\times10^{-06}$\end{tabular} & \begin{tabular}[c]{@{}c@{}}$\pmb{7.58\times10^{-06}}$\end{tabular} &  { \begin{tabular}[c]{@{}c@{}}$8.14\times10^{-06}$\end{tabular}} & \begin{tabular}[c]{@{}l@{}}$\pmb{7.58\times10^{-06}}$\end{tabular} \\ \hline
                        & $E_\theta^{\text{max}}$ & \begin{tabular}[c]{@{}c@{}}$2.49\times10^{-07}$\end{tabular} & \begin{tabular}[c]{@{}c@{}}$2.49\times10^{-07}$\end{tabular} &  \begin{tabular}[c]{@{}c@{}}$\pmb{2.39\times10^{-06}}$\end{tabular} & \begin{tabular}[c]{@{}l@{}}$2.49\times10^{-07}$\end{tabular} \\ \cline{2-6} 
\multirow{-2}{*}{$501$} & $E_\theta^{\text{avg}}$ & \begin{tabular}[c]{@{}c@{}}$1.10\times10^{-08}$\end{tabular} &  \begin{tabular}[c]{@{}c@{}}$\pmb{1.09\times10^{-08}}$\end{tabular} & { \begin{tabular}[c]{@{}c@{}}$1.55\times10^{-08}$\end{tabular}} & \begin{tabular}[c]{@{}l@{}}$\pmb{1.09\times10^{-08}}$\end{tabular} \\ \hline
\end{tabular}
\caption{Comparisons in the case of $\theta=3$ for the function $g(t)=2\exp(-500(t-0.5)^2)+\exp(-3.5t)$ on the interval $[0,1]$.}
\label{Table_theta3_fct2}
\end{table}

\begin{table}[!htbp]
\centering
\begin{tabular}{|Sc|Sc|Sc|Sc|Sc|Sc|}
\cline{1-6}
\multicolumn{2}{|c|}{$\theta=3$} & \multicolumn{4}{c|}{Method}     \\ \hline
$N$                     & Error                   & Method 1                                                     & Method 2                                                                                                                                                                                  & NS                                                                                                      & NAK                                                      \\ \hline
                        & $E_\theta^{\text{max}}$ & \begin{tabular}[c]{@{}c@{}}$1.50\times10^{-01}$\end{tabular} & \begin{tabular}[c]{@{}c@{}}{$4.71\times10^{-02}$}\end{tabular} & { \begin{tabular}[c]{@{}c@{}}$1.13$\end{tabular}}  & \begin{tabular}[c]{@{}l@{}}\pmb{$2.83\times10^{-02}$}\end{tabular} \\ \cline{2-6} 
\multirow{-2}{*}{$31$}  & $E_\theta^{\text{avg}}$ & \begin{tabular}[c]{@{}c@{}}$4.67\times10^{-03}$\end{tabular} & \begin{tabular}[c]{@{}c@{}}$2.61\times10^{-03}$\end{tabular} & { \begin{tabular}[c]{@{}c@{}}$3.42\times10^{-02}$\end{tabular}} & \begin{tabular}[c]{@{}l@{}}\pmb{$1.03\times10^{-03}$}\end{tabular} \\ \hline
                        & $E_\theta^{\text{max}}$ & \begin{tabular}[c]{@{}c@{}}$4.60\times10^{-03}$\end{tabular} & \begin{tabular}[c]{@{}c@{}}$1.31\times10^{-03}$\end{tabular} & { \begin{tabular}[c]{@{}c@{}}$1.07\times10^{-01}$\end{tabular}} & \begin{tabular}[c]{@{}l@{}}\pmb{$2.85\times10^{-04}$}\end{tabular} \\ \cline{2-6} 
\multirow{-2}{*}{$101$} & $E_\theta^{\text{avg}}$ & \begin{tabular}[c]{@{}c@{}}$4.45\times10^{-05}$\end{tabular} & \begin{tabular}[c]{@{}c@{}}$2.39\times10^{-05}$\end{tabular} & { \begin{tabular}[c]{@{}c@{}}$9.95\times10^{-04}$\end{tabular}} & \begin{tabular}[c]{@{}l@{}}\pmb{$4.82\times10^{-06}$}\end{tabular} \\ \hline
                        & $E_\theta^{\text{max}}$ & \begin{tabular}[c]{@{}c@{}}$3.84\times10^{-05}$\end{tabular} & \begin{tabular}[c]{@{}c@{}}$1.06\times10^{-05}$\end{tabular} & { \begin{tabular}[c]{@{}c@{}}$4.36\times10^{-03}$\end{tabular}} & \begin{tabular}[c]{@{}l@{}}\pmb{$4.93\times10^{-07}$}\end{tabular} \\ \cline{2-6} 
\multirow{-2}{*}{$501$} & $E_\theta^{\text{avg}}$ & \begin{tabular}[c]{@{}c@{}}$7.57\times10^{-08}$\end{tabular} & \begin{tabular}[c]{@{}c@{}}$4.00\times10^{-08}$\end{tabular} & { \begin{tabular}[c]{@{}c@{}}$8.14\times10^{-06}$\end{tabular}} & \begin{tabular}[c]{@{}l@{}}\pmb{$5.00\times10^{-09}$}\end{tabular} \\ \hline
\end{tabular}
\caption{Comparisons in the case of $\theta=3$ for the function $g(t)=(t-2)^9+(t-2)^8+(t-2)^4+(t-2)$ on the interval $[0,2]$.}
\label{Table_theta3_fct3}
\end{table}

\begin{table}[!htbp]
\centering
\begin{tabular}{|Sc|Sc|Sc|Sc|Sc|Sc|}
\cline{1-6}
\multicolumn{2}{|c|}{$\theta=3$} & \multicolumn{4}{c|}{Method}     \\ \hline
$N$                     & Error                   & Method 1                                                      & Method 2                                                                                                                                                                                & NS                                                                                                      & NAK                                                      \\ \hline
                        & $E_\theta^{\text{max}}$ & \begin{tabular}[c]{@{}c@{}}$\pmb{1.31\times10^{-03}}$\end{tabular} & \begin{tabular}[c]{@{}c@{}}$\pmb{1.31\times10^{-03}}$\end{tabular} & \begin{tabular}[c]{@{}c@{}}$\pmb{1.31\times10^{-03}}$\end{tabular} & \begin{tabular}[c]{@{}l@{}}$\pmb{1.31\times10^{-03}}$\end{tabular} \\ \cline{2-6} 
\multirow{-2}{*}{$31$}  & $E_\theta^{\text{avg}}$ & \begin{tabular}[c]{@{}c@{}}$6.12\times10^{-05}$\end{tabular} & \begin{tabular}[c]{@{}c@{}}\pmb{$6.10\times10^{-05}$}\end{tabular} & { \begin{tabular}[c]{@{}c@{}}$6.34\times10^{-05}$\end{tabular}} & \begin{tabular}[c]{@{}l@{}}\pmb{$6.10\times10^{-05}$}\end{tabular} \\ \hline
                        & $E_\theta^{\text{max}}$ & \begin{tabular}[c]{@{}c@{}}$\pmb{6.47\times10^{-06}}$\end{tabular} & \begin{tabular}[c]{@{}c@{}}$\pmb{6.47\times10^{-06}}$\end{tabular} & \begin{tabular}[c]{@{}c@{}}$\pmb{6.47\times10^{-06}}$\end{tabular} & \begin{tabular}[c]{@{}l@{}}$\pmb{6.47\times10^{-06}}$\end{tabular} \\ \cline{2-6} 
\multirow{-2}{*}{$101$} & $E_\theta^{\text{avg}}$ & \begin{tabular}[c]{@{}c@{}}$3.23\times10^{-07}$\end{tabular} & \begin{tabular}[c]{@{}c@{}}$3.22\times10^{-07}$\end{tabular} & { \begin{tabular}[c]{@{}c@{}}$3.95\times10^{-07}$\end{tabular}} & \begin{tabular}[c]{@{}l@{}}$\pmb{3.20\times10^{-07}}$\end{tabular} \\ \hline
                        & $E_\theta^{\text{max}}$ & \begin{tabular}[c]{@{}c@{}}$\pmb{9.95\times10^{-09}}$\end{tabular} & \begin{tabular}[c]{@{}c@{}}$\pmb{9.95\times10^{-09}}$\end{tabular} & { \begin{tabular}[c]{@{}c@{}}$1.64\times10^{-07}$\end{tabular}} & \begin{tabular}[c]{@{}l@{}}$\pmb{9.95\times10^{-09}}$\end{tabular} \\ \cline{2-6} 
\multirow{-2}{*}{$501$} & $E_\theta^{\text{avg}}$ & \begin{tabular}[c]{@{}c@{}}$5.05\times10^{-10}$\end{tabular} & \begin{tabular}[c]{@{}c@{}}$5.02\times10^{-10}$\end{tabular} & { \begin{tabular}[c]{@{}c@{}}$1.11\times10^{-09}$\end{tabular}} & \begin{tabular}[c]{@{}l@{}}$\pmb{5.00\times10^{-10}}$\end{tabular} \\ \hline
\end{tabular}
\caption{Comparisons in the case of $\theta=3$ for the function $\displaystyle g(t)=(1+25(t-1)^2)^{-1}$ on the interval $[0,2]$. }
\label{Table_theta3_fct4}
\end{table}

\subsection{Error analysis of quartic splines}
The results obtained for splines of degree $\theta=4$ are illustrated in 
Tables \ref{Table_theta4_fct1}, \ref{Table_theta4_fct2}, \ref{Table_theta4_fct3} and \ref{Table_theta4_fct4}. As these splines are not traditionally used in practice, comparisons with traditional methods are difficult to make. The results of Method 2 will then be compared with the traditional cubic splines with not-a-knot boundary conditions in order to verify how much accuracy is gained by increasing the degree of the spline by one. Let us note however that in some cases, the value of $N$ was too small to get an increase in the overall accuracy of the spline interpolation. This is to be expected when using splines of high degree and a relatively low number of points in order to approximate the high-order derivatives at the nodes. In those cases, not enough information is available for the results to be accurate, so the relative accuracy gain will simply not be computed.

From the results, it can be seen that the accuracy gain increases as $N$ increases. For higher values of $N$ ($N=101$ or $N=501$), the error between the analytic function and the spline is between $6$ to $157$ times smaller when using either Method 1 or Method 2 than with the traditional cubic splines with not-a-knot boundary conditions. While even-degree splines are not usually used in practice, our results suggest that Beaudoin and Beauchemin's spline interpolation method, combined with either of our algorithms (Method 1 or Method 2), may help in cases where accuracy is required at a relatively low computational cost.

\begin{table}[!htbp]
\centering
\begin{tabular}{|Scc|Sc|Sc|Sc|Sc|}
\cline{1-6}
\multicolumn{2}{|c|}{$\theta=4$}             & \multicolumn{4}{c|}{Method}                                       \\ \hline
\multicolumn{1}{|c|}{$N$}                    & Error                  & Method 1            & Method 2            & Cubic spline (NAK)  & Accuracy gain         \\ \hline
\multicolumn{1}{|c|}{\multirow{2}{*}{$31$}}  & $E_\theta^{\text{max}}$ & $1.46\times10^{-03}$ & \pmb{$7.41\times10^{-04}$} & ${3.59\times10^{-03}}$ &$4.84$
\\ \cline{2-6} 
\multicolumn{1}{|c|}{}                       & $E_\theta^{\text{avg}}$ & $5.63\times10^{-05}$ & $\pmb{2.60\times10^{-05}}$ & $1.28\times10^{-04}$ &$4.92$
\\ \hline
\multicolumn{1}{|c|}{\multirow{2}{*}{$101$}} & $E_\theta^{\text{max}}$ & $1.27\times10^{-05}$ & \pmb{$5.55\times10^{-06}$} & ${3.92\times10^{-05}}$ &$7.06$
\\ \cline{2-6} 
\multicolumn{1}{|c|}{}                       & $E_\theta^{\text{avg}}$ & $1.54\times10^{-07}$ & \pmb{$6.46\times10^{-08}$} & ${5.65\times10^{-07}}$ &$8.75$
\\ \hline
\multicolumn{1}{|c|}{\multirow{2}{*}{$501$}} & $E_\theta^{\text{max}}$ & $2.09\times10^{-08}$ & \pmb{$8.96\times10^{-09}$} & ${6.65\times10^{-08}}$ &$7.42$
\\ \cline{2-6} 
\multicolumn{1}{|c|}{}                       & $E_\theta^{\text{avg}}$ & $5.09\times10^{-11}$ & \pmb{$2.14\times10^{-11}$} & ${5.18\times10^{-10}}$ &$24.2$
\\ \hline
\end{tabular}
\caption{Comparisons in the case of $\theta=4$ for the function $g(t)=\sin(3t)e^{-t}$ on the interval $[0,2\pi]$.}
\label{Table_theta4_fct1}
\end{table}

\begin{table}[!htbp]
\centering
\begin{tabular}{|Scc|Sc|Sc|Sc|Sc|}
\cline{1-6}
\multicolumn{2}{|c|}{$\theta=4$}             & \multicolumn{4}{c|}{Method}                                       \\ \hline
\multicolumn{1}{|c|}{$N$}                    & Error                  & Method 1            & Method 2            & Cubic spline (NAK)  & Accuracy gain          \\ \hline
\multicolumn{1}{|c|}{\multirow{2}{*}{$31$}}  & $E_\theta^{\text{max}}$ & $1.92\times10^{-02}$ & $\pmb{1.88\times10^{-02}}$ & $4.29\times10^{-02}$ &$2.28$
\\ \cline{2-6} 
\multicolumn{1}{|c|}{}                       & $E_\theta^{\text{avg}}$ & $4.63\times10^{-03}$ & $5.25\times10^{-03}$ & \pmb{${2.40\times10^{-03}}$} &-
\\ \hline
\multicolumn{1}{|c|}{\multirow{2}{*}{$101$}} & $E_\theta^{\text{max}}$ & $\pmb{1.07\times10^{-05}}$ & $\pmb{1.07\times10^{-05}}$ & ${1.70\times10^{-04}}$ &$15.9$
\\ \cline{2-6} 
\multicolumn{1}{|c|}{}                       & $E_\theta^{\text{avg}}$ & $\pmb{6.89\times10^{-07}}$ & $\pmb{6.89\times10^{-07}}$ & ${7.58\times10^{-06}}$ &$11$
\\ \hline
\multicolumn{1}{|c|}{\multirow{2}{*}{$501$}} & $E_\theta^{\text{max}}$ & $\pmb{2.32\times10^{-09}}$ & $\pmb{2.32\times10^{-09}}$ & ${2.49\times10^{-07}}$ &$107$
\\ \cline{2-6} 
\multicolumn{1}{|c|}{}                       & $E_\theta^{\text{avg}}$ & $\pmb{1.40\times10^{-10}}$ & $\pmb{1.40\times10^{-10}}$ & ${1.09\times10^{-08}}$ &$77.9$
\\ \hline
\end{tabular}
\caption{Comparisons in the case of $\theta=4$ for the function $g(t)=2\exp(-500(t-0.5)^2)+\exp(-3.5t)$ on the interval $[0,1]$.}
\label{Table_theta4_fct2}
\end{table}

\begin{table}[!htbp]
\centering
\begin{tabular}{|Scc|Sc|Sc|Sc|Sc|}
\cline{1-6}
\multicolumn{2}{|c|}{$\theta=4$}             & \multicolumn{4}{c|}{Method}                                       \\ \hline
\multicolumn{1}{|c|}{$N$}                    & Error                  & Method 1            & Method 2            & Cubic spline (NAK)   & Accuracy gain         \\ \hline
\multicolumn{1}{|c|}{\multirow{2}{*}{$31$}}  & $E_\theta^{\text{max}}$ & $1.02\times10^{-02}$ & \pmb{$4.87\times10^{-03}$} & ${2.83\times10^{-02}}$ &$5.81$
\\ \cline{2-6} 
\multicolumn{1}{|c|}{}                       & $E_\theta^{\text{avg}}$ & $3.98\times10^{-04}$ & \pmb{$1.70\times10^{-04}$} & ${1.03\times10^{-03}}$ &$6.06$
\\ \hline
\multicolumn{1}{|c|}{\multirow{2}{*}{$101$}} & $E_\theta^{\text{max}}$ & $9.33\times10^{-05}$ & \pmb{$3.44\times10^{-05}$} & ${2.85\times10^{-04}}$ &$8.28$
\\ \cline{2-6} 
\multicolumn{1}{|c|}{}                       & $E_\theta^{\text{avg}}$ & $1.12\times10^{-06}$ & \pmb{$7.61\times10^{-07}$} & ${4.82\times10^{-06}}$ &$6.33$
\\ \hline
\multicolumn{1}{|c|}{\multirow{2}{*}{$501$}} & $E_\theta^{\text{max}}$ & $1.55\times10^{-07}$ & \pmb{$6.71\times10^{-08}$} & ${4.93\times10^{-07}}$ &$7.35$
\\ \cline{2-6} 
\multicolumn{1}{|c|}{}                       & $E_\theta^{\text{avg}}$ & $3.77\times10^{-10}$ & \pmb{$1.56\times10^{-10}$} & ${5.00\times10^{-09}}$ &$32.1$
\\ \hline
\end{tabular}
\caption{Comparisons in the case of $\theta=4$ for the function $g(t)=(t-2)^9+(t-2)^8+(t-2)^4+(t-2)$ on the interval $[0,2]$. }
\label{Table_theta4_fct3}
\end{table}

\begin{table}[!htbp]
\centering
\begin{tabular}{|Scc|Sc|Sc|Sc|Sc|}
\cline{1-6}
\multicolumn{2}{|c|}{$\theta=4$}             & \multicolumn{4}{c|}{Method}                                       \\ \hline
\multicolumn{1}{|c|}{$N$}                    & Error                  & Method 1            & Method 2            & Cubic spline (NAK)  & Accuracy gain          \\ \hline
\multicolumn{1}{|c|}{\multirow{2}{*}{$31$}}  & $E_\theta^{\text{max}}$ & $4.59\times10^{-04}$ & \pmb{$4.53\times10^{-04}$} & ${1.31\times10^{-03}}$ &$2.89$
\\ \cline{2-6} 
\multicolumn{1}{|c|}{}                       & $E_\theta^{\text{avg}}$ & $7.37\times10^{-05}$ & $8.26\times10^{-05}$ & \pmb{${6.10\times10^{-05}}$} & -
\\ \hline
\multicolumn{1}{|c|}{\multirow{2}{*}{$101$}} & $E_\theta^{\text{max}}$ & $\pmb{2.66\times10^{-07}}$ & $2.67\times10^{-07}$ & $6.47\times10^{-06}$ & $24.2$
\\ \cline{2-6} 
\multicolumn{1}{|c|}{}                       & $E_\theta^{\text{avg}}$ & $\pmb{1.74\times10^{-08}}$ & $\pmb{1.74\times10^{-08}}$ & ${3.20\times10^{-07}}$ &$18.4$
\\ \hline
\multicolumn{1}{|c|}{\multirow{2}{*}{$501$}} & $E_\theta^{\text{max}}$ & $\pmb{6.33\times10^{-11}}$ & $\pmb{6.33\times10^{-11}}$ & $9.95\times10^{-09}$ &$157$
\\ \cline{2-6} 
\multicolumn{1}{|c|}{}                       & $E_\theta^{\text{avg}}$ & $4.18\times10^{-12}$ & $\pmb{4.16\times10^{-12}}$ & $5.00\times10^{-10}$ &$120$
\\ \hline
\end{tabular}
\caption{Comparisons in the case of $\theta=4$ for the function $\displaystyle g(t)=(1+25(t-1)^2)^{-1}$ on the interval $[0,2]$. }
\label{Table_theta4_fct4}
\end{table}

\subsection{Error analysis of quintic splines}
Tables \ref{Table_theta5_fct1}, \ref{Table_theta5_fct2}, \ref{Table_theta5_fct3} and \ref{Table_theta5_fct4} illustrate the results obtained when using splines of degree $\theta=5$. For the same reasons as $\theta=4$, we are not able to compare this method with traditional methods of the same degree. Comparisons with the traditional cubic spline with a not-a-knot boundary condition will however be made. 

The results show that, again, the gain in accuracy increases as $N$ increases. However, in this case, the error between the analytic function and the spline function is between $65$ to $10\,000$ times smaller than the results obtained using the traditional cubic spline (with a not-a-knot boundary conditions) for higher values of $N$ ($N=101$ or $N=501$). This represents a significant gain in accuracy by simply increasing the degree of the spline by two.

\begin{table}[!htbp]
\centering
\begin{tabular}{|Scc|Sc|Sc|Sc|Sc|}
\cline{1-6}
\multicolumn{2}{|c|}{$\theta=5$}             & \multicolumn{4}{c|}{Method}                \\ \hline
\multicolumn{1}{|c|}{$N$}                    & Error                  & Method 1            & Method 2     & Cubic spline (NAK) & Accuracy gain      \\ \hline
\multicolumn{1}{|c|}{\multirow{2}{*}{$31$}}  & $E_\theta^{\text{max}}$ & $5.52\times10^{-04}$ & $\pmb{6.24\times10^{-05}}$ & ${3.59\times10^{-03}}$ &$57.5$
\\ \cline{2-6} 
\multicolumn{1}{|c|}{}                       & $E_\theta^{\text{avg}}$ & $1.80\times10^{-05}$ & $\pmb{2.16\times10^{-06}}$ & $1.28\times10^{-04}$ &$59.3$
\\ \hline
\multicolumn{1}{|c|}{\multirow{2}{*}{$101$}} & $E_\theta^{\text{max}}$ & $6.08\times10^{-07}$ & $\pmb{6.30\times10^{-08}}$ & ${3.92\times10^{-05}}$ &$622$
\\ \cline{2-6} 
\multicolumn{1}{|c|}{}                       & $E_\theta^{\text{avg}}$ & $6.20\times10^{-09}$ & $\pmb{8.66\times10^{-10}}$ & ${5.65\times10^{-07}}$ &$652$
\\ \hline
\multicolumn{1}{|c|}{\multirow{2}{*}{$501$}} & $E_\theta^{\text{max}}$ & $7.15\times10^{-11}$ & $\pmb{7.06\times10^{-12}}$ & ${6.65\times10^{-08}}$ &$9.42\times10^{03}$
\\ \cline{2-6} 
\multicolumn{1}{|c|}{}                       & $E_\theta^{\text{avg}}$ & $1.54\times10^{-13}$ & $\pmb{2.72\times10^{-14}}$ & ${5.18\times10^{-10}}$ &$1.00\times10^{04}$
\\ \hline
\end{tabular}
\caption{Comparisons in the case of $\theta=5$ for the function $g(t)=\sin(3t)e^{-t}$ on the interval $[0,2\pi]$.}
\label{Table_theta5_fct1}
\end{table}

\begin{table}[!htbp]
\centering
\begin{tabular}{|Scc|Sc|Sc|Sc|Sc|}
\cline{1-6}
\multicolumn{2}{|c|}{$\theta=5$}             & \multicolumn{4}{c|}{Method}                \\ \hline
\multicolumn{1}{|c|}{$N$}                    & Error                  & Method 1            & Method 2          & Cubic spline (NAK) & Accuracy gain \\ \hline
\multicolumn{1}{|c|}{\multirow{2}{*}{$31$}}  & $E_\theta^{\text{max}}$ & $\pmb{2.09\times10^{-02}}$ & $\pmb{2.09\times10^{-02}}$ & $4.29\times10^{-02}$ &$2.05$
\\ \cline{2-6} 
\multicolumn{1}{|c|}{}                       & $E_\theta^{\text{avg}}$ & $\pmb{1.94\times10^{-03}}$ & $2.68\times10^{-03}$ & ${2.40\times10^{-03}}$ &-
\\ \hline
\multicolumn{1}{|c|}{\multirow{2}{*}{$101$}} & $E_\theta^{\text{max}}$ & $\pmb{2.61\times10^{-06}}$ & $\pmb{2.61\times10^{-06}}$ & ${1.70\times10^{-04}}$ &$65.1$
\\ \cline{2-6} 
\multicolumn{1}{|c|}{}                       & $E_\theta^{\text{avg}}$ & $1\pmb{.15\times10^{-07}}$ & $\pmb{1.15\times10^{-07}}$ & ${7.58\times10^{-06}}$ &$65.9$
\\ \hline
\multicolumn{1}{|c|}{\multirow{2}{*}{$501$}} & $E_\theta^{\text{max}}$ & $\pmb{1.25\times10^{-10}}$ & $\pmb{1.25\times10^{-10}}$ & ${2.49\times10^{-07}}$ &$1.99\times10^{03}$
\\ \cline{2-6} 
\multicolumn{1}{|c|}{}                       & $E_\theta^{\text{avg}}$ & $\pmb{5.20\times10^{-12}}$ & $\pmb{5.20\times10^{-12}}$ & ${1.09\times10^{-08}}$ &$2.10\times10^{03}$
\\ \hline
\end{tabular}
\caption{Comparisons in the case of $\theta=5$ for the function $g(t)=2\exp(-500(t-0.5)^2)+\exp(-3.5t)$ on the interval $[0,1]$. }
\label{Table_theta5_fct2}
\end{table}

\begin{table}[!htbp]
\centering
\begin{tabular}{|Scc|Sc|Sc|Sc|Sc|}
\cline{1-6}
\multicolumn{2}{|c|}{$\theta=5$}             & \multicolumn{4}{c|}{Method}               \\ \hline
\multicolumn{1}{|c|}{$N$}                    & Error                  & Method 1            & Method 2      & Cubic spline (NAK)      & Accuracy gain\\ \hline
\multicolumn{1}{|c|}{\multirow{2}{*}{$31$}}  & $E_\theta^{\text{max}}$ & $1.80\times10^{-03}$ & $\pmb{3.17\times10^{-04}}$ & ${2.83\times10^{-02}}$  &$8.93$
\\ \cline{2-6} 
\multicolumn{1}{|c|}{}                       & $E_\theta^{\text{avg}}$ & $5.57\times10^{-05}$ & $\pmb{1.32\times10^{-05}}$ & ${1.03\times10^{-03}}$ &$78.0$
\\ \hline
\multicolumn{1}{|c|}{\multirow{2}{*}{$101$}} & $E_\theta^{\text{max}}$ & $5.40\times10^{-06}$ & $\pmb{7.98\times10^{-07}}$ & ${2.85\times10^{-04}}$ &$357$
\\ \cline{2-6} 
\multicolumn{1}{|c|}{}                       & $E_\theta^{\text{avg}}$ & $5.17\times10^{-08}$ & $\pmb{1.08\times10^{-08}}$ & ${4.82\times10^{-06}}$  &$446$
\\ \hline
\multicolumn{1}{|c|}{\multirow{2}{*}{$501$}} & $E_\theta^{\text{max}}$ & $1.86\times10^{-09}$ & $\pmb{2.58\times10^{-10}}$ & ${4.93\times10^{-07}}$ &$1.91\times10^{03}$
\\ \cline{2-6} 
\multicolumn{1}{|c|}{}                       & $E_\theta^{\text{avg}}$ & $3.73\times10^{-12}$ & $\pmb{7.18\times10^{-13}}$ & ${5.00\times10^{-09}}$ &$6.96\times10^{03}$
\\ \hline
\end{tabular}
\caption{Comparisons in the case of $\theta=5$ for the function $g(t)=(t-2)^9+(t-2)^8+(t-2)^4+(t-2)$ on the interval $[0,2]$. }
\label{Table_theta5_fct3}
\end{table}

\begin{table}[!htbp]
\centering
\begin{tabular}{|Scc|Sc|Sc|Sc|Sc|}
\cline{1-6}
\multicolumn{2}{|c|}{$\theta=5$}             & \multicolumn{4}{c|}{Method}                \\ \hline
\multicolumn{1}{|c|}{$N$}                    & Error                  & Method 1            & Method 2     & Cubic spline (NAK)   & Accuracy gain    \\ \hline
\multicolumn{1}{|c|}{\multirow{2}{*}{$31$}}  & $E_\theta^{\text{max}}$ & $\pmb{3.75\times10^{-04}}$ & $6.00\times10^{-04}$ & ${1.31\times10^{-03}}$ &$2.18$
\\ \cline{2-6} 
\multicolumn{1}{|c|}{}                       & $E_\theta^{\text{avg}}$ & $\pmb{2.60\times10^{-05}}$ & $6.50\times10^{-05}$ & ${6.10\times10^{-05}}$ &-
\\ \hline
\multicolumn{1}{|c|}{\multirow{2}{*}{$101$}} & $E_\theta^{\text{max}}$ & $\pmb{5.90\times10^{-08}}$ & $\pmb{5.90\times10^{-08}}$ & $6.47\times10^{-06}$  &$110$
\\ \cline{2-6} 
\multicolumn{1}{|c|}{}                       & $E_\theta^{\text{avg}}$ & $2.26\times10^{-09}$ & $\pmb{2.25\times10^{-09}}$ & ${3.20\times10^{-07}}$  &$142$
\\ \hline
\multicolumn{1}{|c|}{\multirow{2}{*}{$501$}} & $E_\theta^{\text{max}}$ & $\pmb{2.98\times10^{-12}}$ & $3.00\times10^{-12}$ & $9.95\times10^{-09}$ &$3.32\times10^{03}$
\\ \cline{2-6} 
\multicolumn{1}{|c|}{}                       & $E_\theta^{\text{avg}}$ & $1.22\times10^{-13}$ & $\pmb{1.15\times10^{-13}}$ & $5.00\times10^{-10}$ &$4.35\times10^{03}$
\\ \hline
\end{tabular}
\caption{Comparisons in the case of $\theta=5$ for the function $\displaystyle g(t)=(1+25(t-1)^2)^{-1}$ on the interval $[0,2]$. }
\label{Table_theta5_fct4}
\end{table}

\subsection{Error analysis of higher degree splines}
Using traditional methods, high-degree spline interpolation is not an easy task due to the lack of effective methods to approximate high-order derivatives of an unknown function at the interpolation nodes. The advantage of the approach presented in this paper is that very accurate interpolation results can be obtained by increasing the degree of the spline. This is clearly an asset for problems where higher accuracy is essential.

Tables \ref{Table_O11_function1}, \ref{Table_O11_function2} and \ref{Table_O11_function3} show the results obtained using $\theta=11$. No results are shown for the analytic function $g(t)=(t-2)^9+(t-2)^8+(t-2)^4+(t-2)$ as the approximation is exact in this case.  This is due to the fact that the function considered is a polynomial of degree nine and we are using polynomials of degree eleven to approximate it. By looking at the results, the same trends can be observed; the accuracy gain increases as $N$ increases. For higher values of $N$ ($N=101$ or $N=501$), the error between the analytic function and the spline function is now between $1.75\times10^{4}$ to $1.89\times10^{13}$ times smaller than the results obtained using the traditional cubic spline. This again represents a significant gain in accuracy. 

\begin{table}[!htbp]
\centering
\begin{tabular}{|Sc|Sc|Sc|Sc|Sc|Sc|}
\hline
\multicolumn{2}{|c|}{$\theta=11$}                & \multicolumn{4}{c|}{Method}          \\ \hline
$N$                    & Error                  & Method 1            & Method 2     & Cubic spline (NAK)       & Accuracy gain \\ \hline
\multirow{2}{*}{$31$}  & $E_\theta^{\text{max}}$ & $4.11\times10^{-06}$ & $\pmb{3.87\times10^{-07}}$ & ${3.59\times10^{-03}}$ &$9.28\times10^{03}$
\\ \cline{2-6} 
                       & $E_\theta^{\text{avg}}$ & $9.37\times10^{-08}$ & $\pmb{9.01\times10^{-09}}$ & $1.28\times10^{-04}$  &$1.42\times10^{04}$
\\ \hline
\multirow{2}{*}{$101$} & $E_\theta^{\text{max}}$ & $1.13\times10^{-11}$ & $\pmb{7.15\times10^{-13}}$ & ${3.92\times10^{-05}}$ &$5.48\times10^{07}$
\\ \cline{2-6} 
                       & $E_\theta^{\text{avg}}$ & $7.80\times10^{-14}$  & $\pmb{6.21\times10^{-15}}$ & ${5.65\times10^{-07}}$  &$9.10\times10^{07}$
\\ \hline
\multirow{2}{*}{$501$} & $E_\theta^{\text{max}}$ & $4.67\times10^{-19}$ & $\pmb{1.34\times10^{-20}}$ & ${6.65\times10^{-08}}$ &$4.96\times10^{12}$
\\ \cline{2-6} 
                       & $E_\theta^{\text{avg}}$ & $6.54\times10^{-22}$ & $\pmb{2.74\times10^{-23}}$ & ${5.18\times10^{-10}}$ &$1.89\times10^{13}$
\\ \hline
\end{tabular}
\caption{Comparisons in the case of $\theta=11$ for the function $g(t)=\sin(3t)e^{-t}$ on the interval $[0,2\pi]$. }
\label{Table_O11_function1}
\end{table}

\begin{table}[!htbp]
\centering
\begin{tabular}{|Sc|Sc|Sc|Sc|Sc|Sc|}
\hline
\multicolumn{2}{|c|}{$\theta=11$}                & \multicolumn{4}{c|}{Method}          \\ \hline
$N$                    & Error                  & Method 1            & Method 2    & Cubic spline (NAK)       & Accuracy gain \\ \hline
\multirow{2}{*}{$31$}  & $E_\theta^{\text{max}}$ & ${7.43\times10^{-01}}$ & $32.2$               & \pmb{$4.29\times10^{-02}$} & -
\\ \cline{2-6} 
                       & $E_\theta^{\text{avg}}$ & ${3.73\times10^{-02}}$ & $1.42$               & \pmb{${2.40\times10^{-03}}$} & -
\\ \hline
\multirow{2}{*}{$101$} & $E_\theta^{\text{max}}$ & $\pmb{1.17\times10^{-10}}$ & $\pmb{1.17\times10^{-10}}$ & ${1.70\times10^{-04}}$ &$1.45\times10^{06}$
\\ \cline{2-6} 
                       & $E_\theta^{\text{avg}}$ & $\pmb{5.51\times10^{-12}}$ & $\pmb{5.51\times10^{-12}}$ & ${7.58\times10^{-06}}$ &$1.38\times10^{06}$
\\ \hline
\multirow{2}{*}{$501$} & $E_\theta^{\text{max}}$ & $\pmb{9.71\times10^{-20}}$ & $\pmb{9.71\times10^{-20}}$ & ${2.49\times10^{-07}}$ &$2.56\times10^{12}$
\\ \cline{2-6} 
                       & $E_\theta^{\text{avg}}$ & $\pmb{4.08\times10^{-21}}$ & $\pmb{4.08\times10^{-21}}$ & ${1.09\times10^{-08}}$ &$2.67\times10^{12}$
\\ \hline
\end{tabular}
\caption{Comparisons in the case of $\theta=11$ for the function $g(t)=2\exp(-500(t-0.5)^2)+\exp(-3.5t)$ on the interval $[0,1]$. }
\label{Table_O11_function2}
\end{table}

\begin{table}[!htbp]
\centering
\begin{tabular}{|Sc|Sc|Sc|Sc|Sc|Sc|}
\hline
\multicolumn{2}{|c|}{$\theta=11$}                & \multicolumn{4}{c|}{Method}          \\ \hline
$N$                    & Error                  & Method 1            & Method 2   & Cubic spline (NAK)   & Accuracy gain      \\ \hline
\multirow{2}{*}{$31$}  & $E_\theta^{\text{max}}$ & ${9.84\times10^{-03}}$ & $4.71\times10^{-01}$ & \pmb{${1.31\times10^{-03}}$} & -
\\ \cline{2-6} 
                       & $E_\theta^{\text{avg}}$ & ${4.94\times10^{-04}}$ & $2.08\times10^{-02}$ & \pmb{${6.10\times10^{-05}}$} & -
\\ \hline
\multirow{2}{*}{$101$} & $E_\theta^{\text{max}}$ & $\pmb{6.08\times10^{-12}}$ & $3.70\times10^{-10}$  & $6.47\times10^{-06}$ &$1.75\times10^{04}$
\\ \cline{2-6} 
                       & $E_\theta^{\text{avg}}$ & $\pmb{1.37\times10^{-13}}$ & $5.13\times10^{-12}$ & ${3.20\times10^{-07}}$ &$6.24\times10^{04}$
\\ \hline
\multirow{2}{*}{$501$} & $E_\theta^{\text{max}}$ & \pmb{${2.33\times10^{-21}}$} & $\pmb{2.33\times10^{-21}}$ & $9.95\times10^{-09}$ &$4.27\times10^{12}$
\\ \cline{2-6} 
                       & $E_\theta^{\text{avg}}$ & $6.51\times10^{-23}$ & $\pmb{6.31\times10^{-23}}$ & $5.00\times10^{-10}$ &$7.92\times10^{12}$
\\ \hline
\end{tabular}
\caption{Comparisons in the case of $\theta=11$ for the function $\displaystyle g(t)=(1+25(t-1)^2)^{-1}$ on the interval $[0,2]$. }
\label{Table_O11_function3}
\end{table}

\subsection{General remarks}
General remarks can be made from the results shown in Tables \ref{Table_theta3_fct1} to \ref{Table_O11_function3}.

\begin{enumerate}
\item[$\bullet$] Problems can arise when using high values of $\theta$ if the number of interpolation nodes used is too small. This is due to the fact that the algorithm cannot accurately approximate the boundary conditions of high-order derivatives if limited information on the function is known.
Users should therefore make sure that an appropriate number of interpolation points are available when using greater values of $\theta$.

\item[$\bullet$] For a chosen value of $\theta$, increasing $\theta$ increases the overall accuracy of the interpolation for both proposed methods. 

\item[$\bullet$] As $\theta$ increases, the number of digits used for the computation should be increased in order to be able to adequately measure the accuracy obtained.

\item[$\bullet$] The time-complexity of the algorithms increases as $\theta$ increases. However, for problems where accuracy gain is more important than time-complexity, our two algorithms have proved to yield very accurate results.

\item[$\bullet$] Generally, the gain in accuracy for both algorithms presented (Method 1 and Method 2) is very similar. For the sake of simplicity and because it can cope with any integer value of $N$, Method 1 may be preferred over Method 2.

\item[$\bullet$] While traditional cubic splines with natural or not-a-knot boundary conditions may offer satisfying results in many cases, a significant gain in accuracy can be obtained by simply using our proposed method with a higher value of $\theta$. This can be achieved without the need of additional information on the function to interpolate. 

\item[$\bullet$] While high-degree splines (such as $\theta=11$) leads to very accurate results, splines of degree 4 or 5 may be a good compromise if the user only has limited information available.

\item[$\bullet$] Increasing $N$ leads to an exponential growth in gain of accuracy when comparing cubic spline interpolation with not-a-knot boundary conditions and splines of degree $\theta=11$. The relative accuracy gain in terms of $N$ for functions $g(t)=2\exp(-500(t-0.5)^2)+\exp(-3.5t)$ on the interval $[0,1]$  and $g(t)=(1+25(t-1)^2)^{-1}$ on the interval $[0,2]$ are illustrated in Figures \ref{fig:AccGain_ex1} and \ref{fig:AccGain_ex2} respectively.

\begin{figure}[!htbp]
  \centering
  \includegraphics[width=0.5\linewidth,scale=0.5]{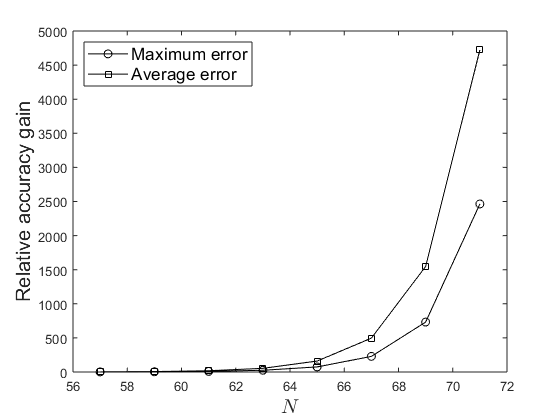}
  \caption{Relative accuracy gain in terms of $N$ for the function $g(t)=2\exp(-500(t-0.5)^2)+\exp(-3.5t)$ on the interval $[0,1]$.}
  \label{fig:AccGain_ex1}
\end{figure}

\begin{figure}[!htbp]
  \centering
  \includegraphics[width=0.5\linewidth,scale=0.5]{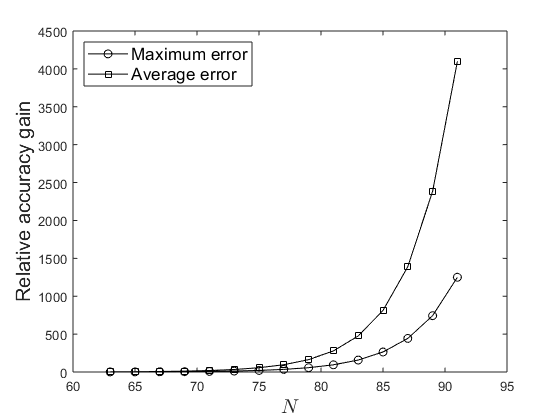}
  \caption{Relative accuracy gain in terms of $N$ for the function $g(t)=(1+25(t-1)^2)^{-1}$ on the interval $[0,2]$.}
  \label{fig:AccGain_ex2}
\end{figure}

\end{enumerate}

\section{Practical example}
\label{Sec:Practical_example}
In \cite{beaudoin_tutorial/article_1998}, the authors have concluded that their method leads to an accurate approximation of the analytical Fourier transform of a digitized function. The goal of this practical example is to pair this method with the proposed algorithm to compute the boundary conditions (Method 1). In the following, we consider $N=2^{13}$ points of the digitized function $f(t)=\cos(60t)e^{-2t}$ on the interval $[0,81.92]$. To show that accurate boundary conditions are essential, we calculate an approximation of the Fourier transform of an analytical function (known for comparison purposes) from a digitized version of this function with the method used in \cite{beaudoin_tutorial/article_1998}, which is essentially based on spline interpolation of different degrees. The effect of using more accurate boundary conditions (e.g. Method 1) combined with the effect of increasing the degree of the spline interpolation on the accuracy of the approximation of the Fourier transform is shown and the numerical results are compared with the analytical results. Let us note that the FFT algorithm used in the process can be found in \cite{brigham_fast_1988}.

\begin{figure}
    \centering
    \includegraphics[width=1\linewidth,scale=0.5]{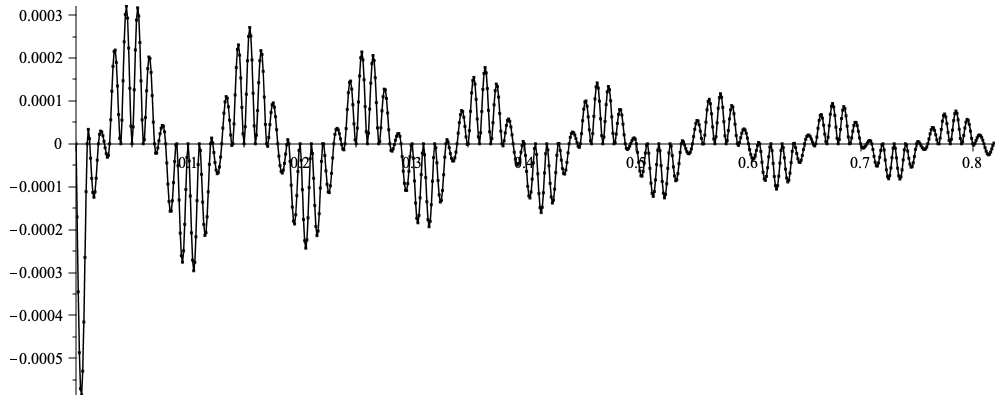}
    \caption{Difference between the function $f(t)$ and its cubic spline interpolation with exact boundary conditions.}
    \label{figA}
\end{figure}

\begin{figure}
    \centering
    \includegraphics[width=1\linewidth,scale=0.5]{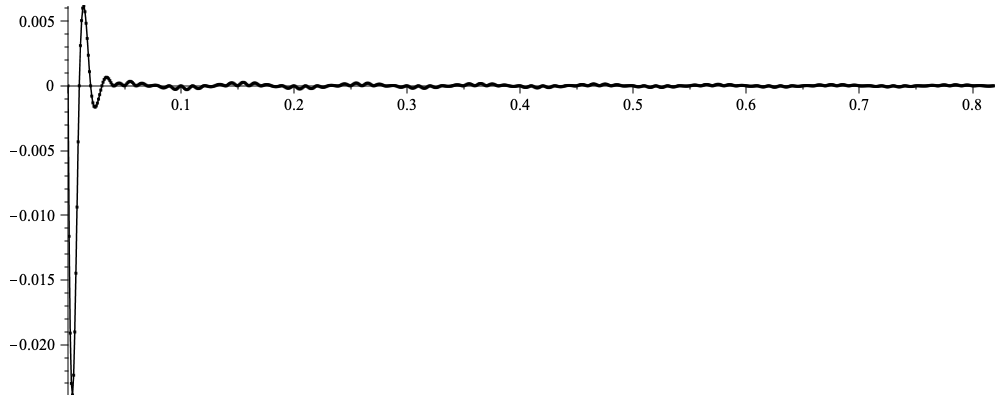}
    \caption{Difference between the function $f(t)$ and its cubic spline interpolation with boundary conditions $b_1=b_2=0$.}
    \label{figB}
\end{figure}

In Figures \ref{figA} and \ref{figB}, where cubic splines were used, the difference between the analytical function and the spline interpolation of the digitized function shows that the interpolation is not curved enough between the sampling points. This produces a  spurious peak in the frequency domain at the sampling frequency and at the frequency of the function itself, which then leads to an error on the computed Fourier transform. If the boundary conditions are not accurate enough, this error is amplified.

\begin{figure}[htp]

\begin{minipage}{\linewidth}
  \centering
  \includegraphics[width=1\textwidth]{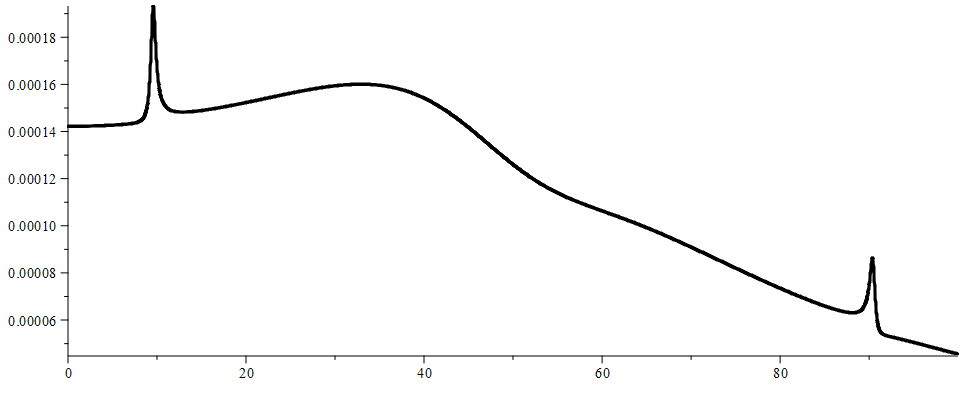}
\end{minipage}
\vspace{1em} 

\begin{minipage}{\linewidth}
  \centering
  \includegraphics[width=1\textwidth]{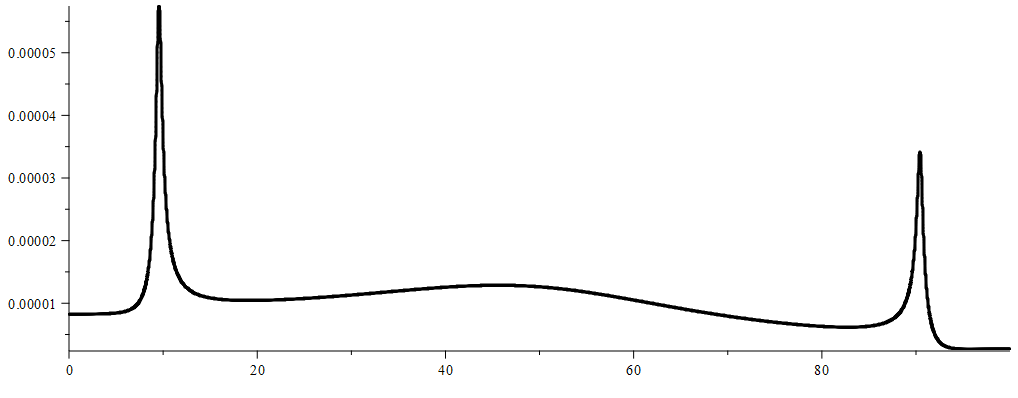}
\end{minipage}
\vspace{1em} 

\begin{minipage}{\linewidth}
  \centering
  \includegraphics[width=1\textwidth]{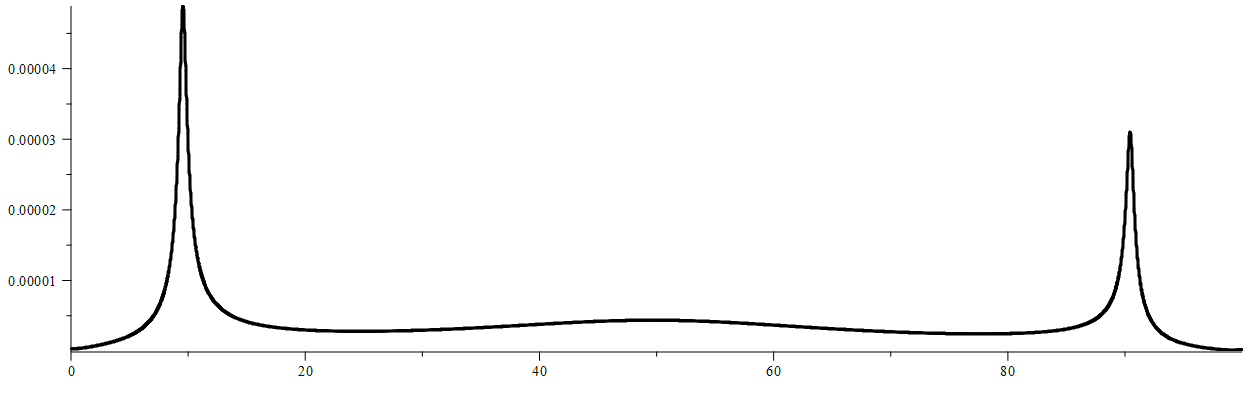}
\end{minipage}
\vspace{1em} 

\caption{Absolute difference in terms of $t$ between the approximation of the analytical Fourier transform (using splines of degree 3) and the exact analytical Fourier transform of the function $f(t)$. First figure : Boundary conditions were set $b_n=0$, $n=1,2$, and $b_0=f(81.92)-f(0)$. Second figure : Boundary conditions were computed with Method 1. Third figure : Boundary conditions were computed analytically.}
\label{fig:figuredeg3}

\end{figure}

\begin{figure}[htp]

\begin{minipage}{\linewidth}
  \centering
  \includegraphics[width=1\textwidth]{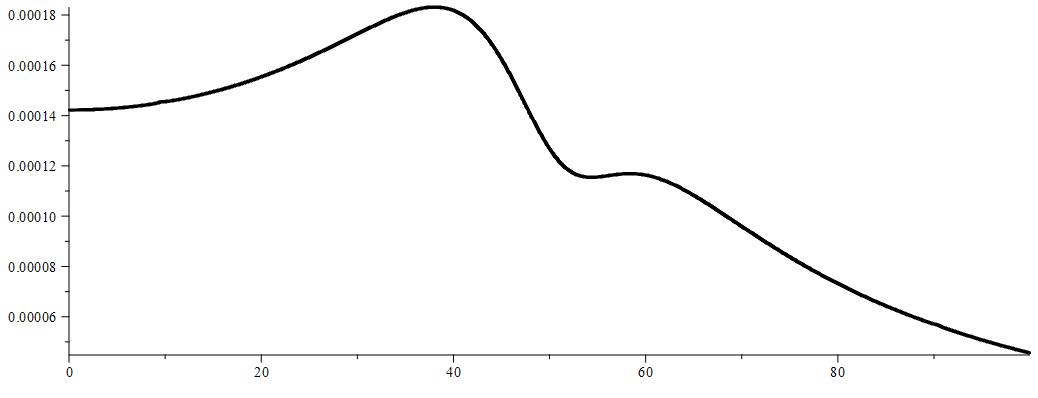}
\end{minipage}
\vspace{1em} 

\begin{minipage}{\linewidth}
  \centering
  \includegraphics[width=1\textwidth]{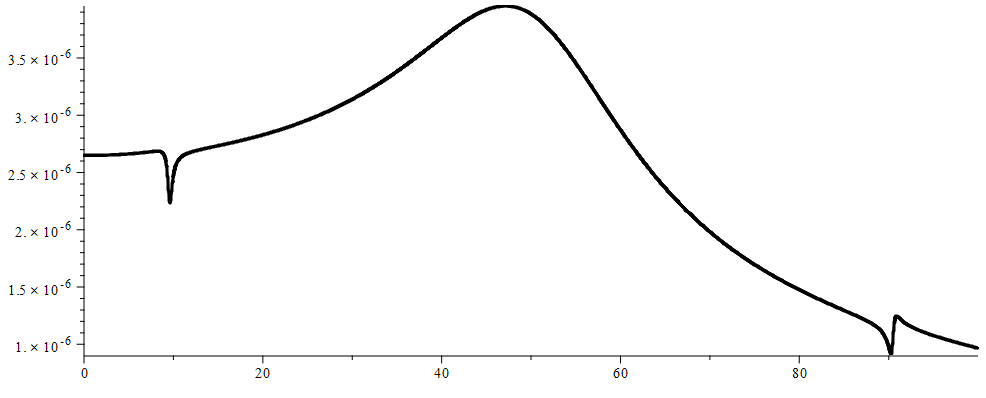}
\end{minipage}
\vspace{1em} 

\begin{minipage}{\linewidth}
  \centering
  \includegraphics[width=1\textwidth]{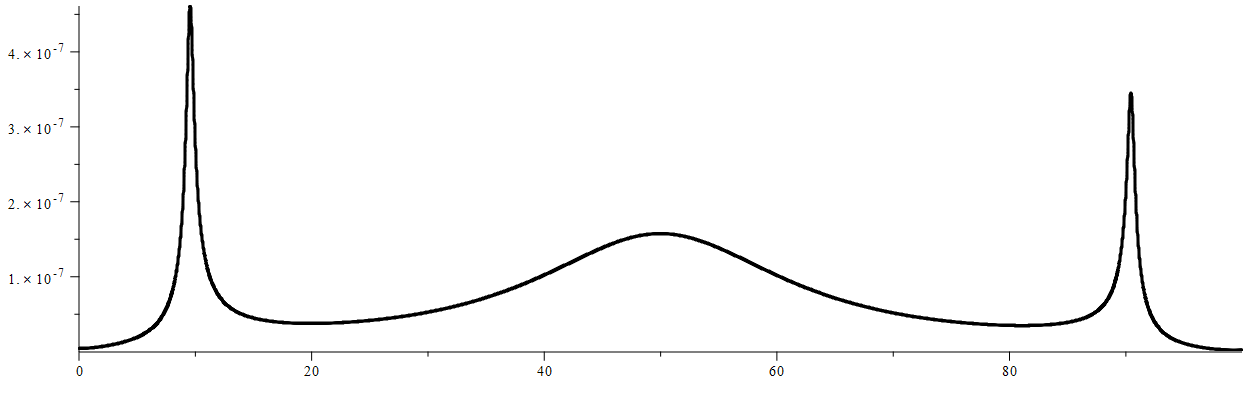}
\end{minipage}
\vspace{1em} 

\caption{Absolute difference in terms of $t$ between the approximation of the analytical Fourier transform (using splines of degree 5) and the exact analytical Fourier transform of the function $f(t)$. First figure : Boundary conditions were set $b_n=0$, $n=1,\ldots,4$, and $b_0=f(81.92)-f(0)$. Second figure : Boundary conditions were computed with Method 1. Third figure : Boundary conditions were computed analytically.}
\label{fig:figuredeg5}

\end{figure}

\begin{figure}[htp]

\begin{minipage}{\linewidth}
  \centering
  \includegraphics[width=1\textwidth]{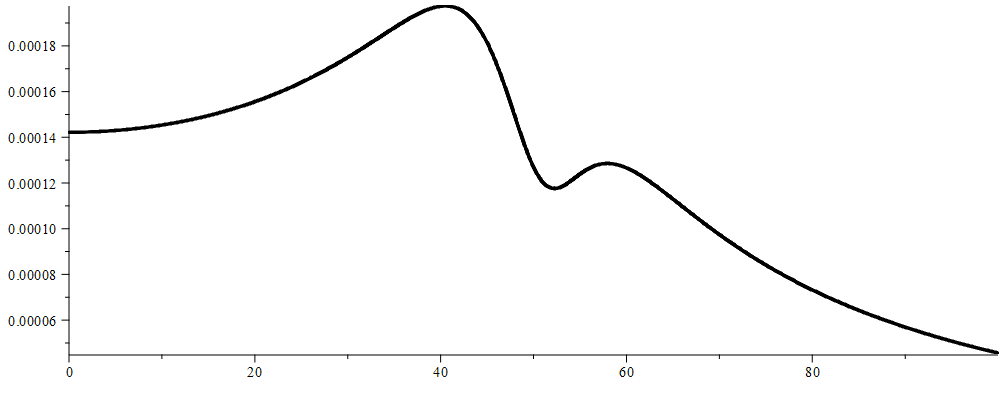}
\end{minipage}
\vspace{1em} 

\begin{minipage}{\linewidth}
  \centering
  \includegraphics[width=1\textwidth]{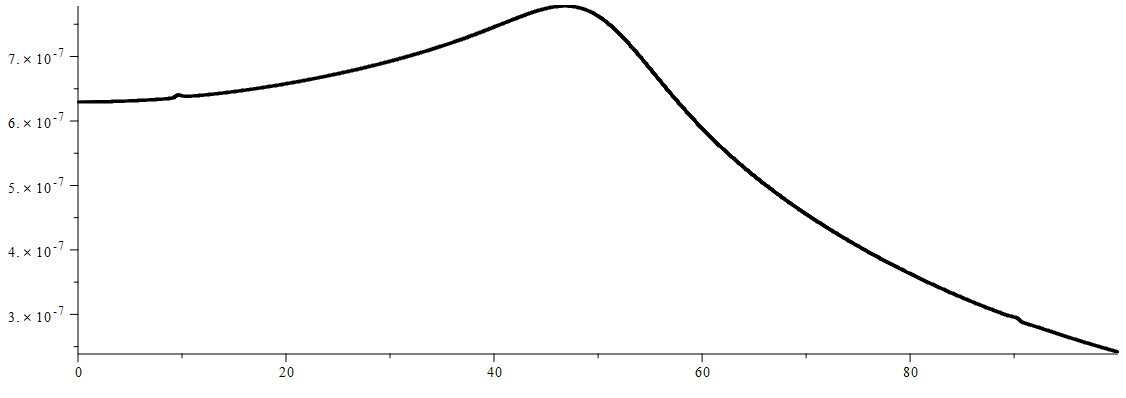}
\end{minipage}
\vspace{1em} 

\begin{minipage}{\linewidth}
  \centering
  \includegraphics[width=1\textwidth]{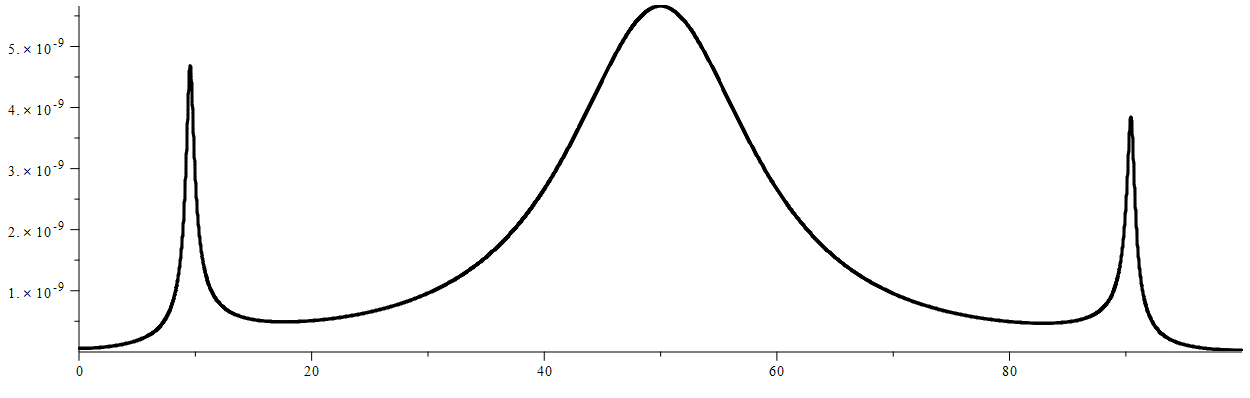}
\end{minipage}
\vspace{1em} 

\caption{Absolute difference in terms of $t$ between the approximation of the analytical Fourier transform (using splines of degree 7) and the exact analytical Fourier transform of the function $f(t)$. First figure : Boundary conditions were set $b_n=0$, $n=1,\ldots,6$, and $b_0=f(81.92)-f(0)$. Second figure : Boundary conditions were computed with Method 1. Third figure : Boundary conditions were computed analytically.}
\label{fig:figuredeg7}

\end{figure}

\begin{figure}[htp]

\begin{minipage}{\linewidth}
  \centering
  \includegraphics[width=1\textwidth]{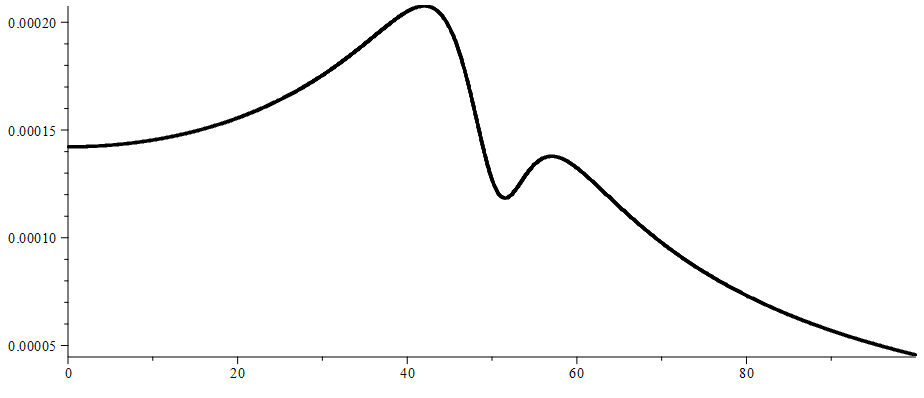}
\end{minipage}
\vspace{1em} 

\begin{minipage}{\linewidth}
  \centering
  \includegraphics[width=1\textwidth]{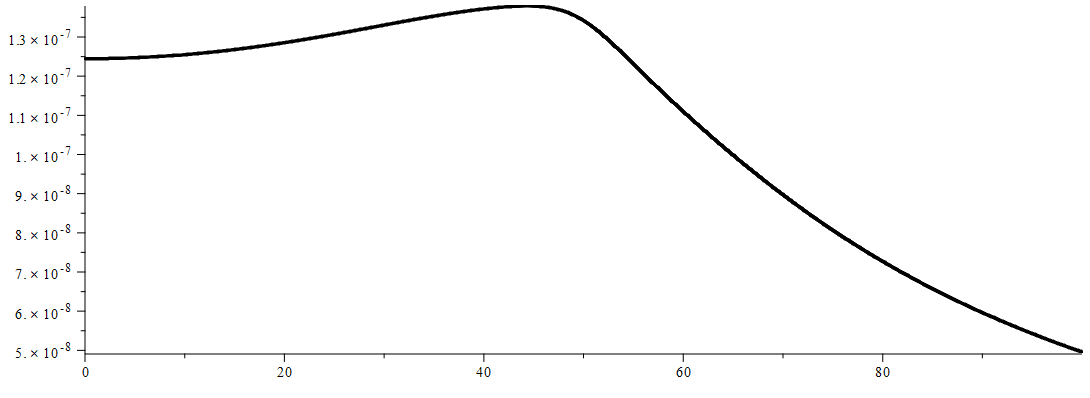}
\end{minipage}
\vspace{1em} 

\begin{minipage}{\linewidth}
  \centering
  \includegraphics[width=1\textwidth]{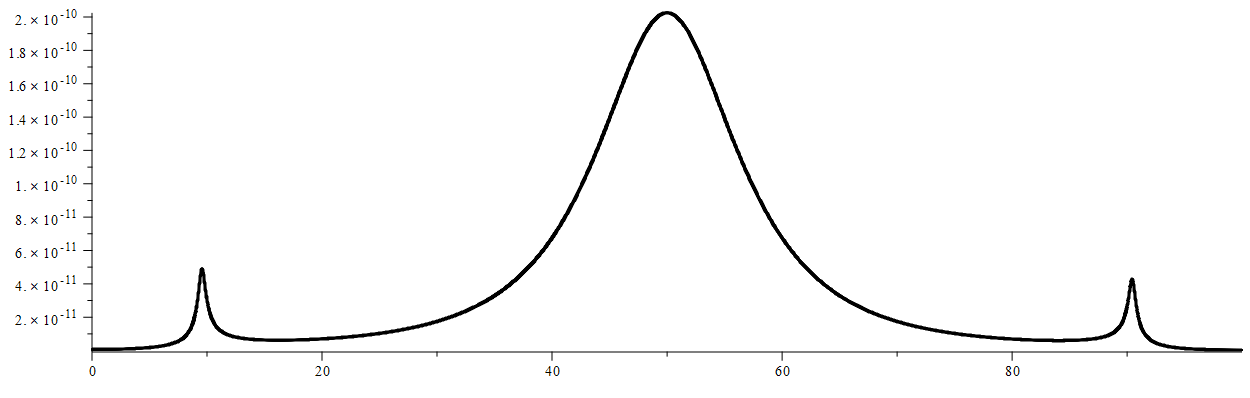}
\end{minipage}
\vspace{1em} 

\caption{Absolute difference in terms of $t$ between the approximation of the analytical Fourier transform (using splines of degree 9) and the exact analytical Fourier transform of the function $f(t)$. First figure : Boundary conditions were set $b_n=0$, $n=1,\ldots,8$, and $b_0=f(81.92)-f(0)$. Second figure : Boundary conditions were computed with Method 1. Third figure : Boundary conditions were computed analytically.}
\label{fig:figuredeg9}

\end{figure}

\begin{figure}[htp]

\begin{minipage}{\linewidth}
  \centering
  \includegraphics[width=1\textwidth]{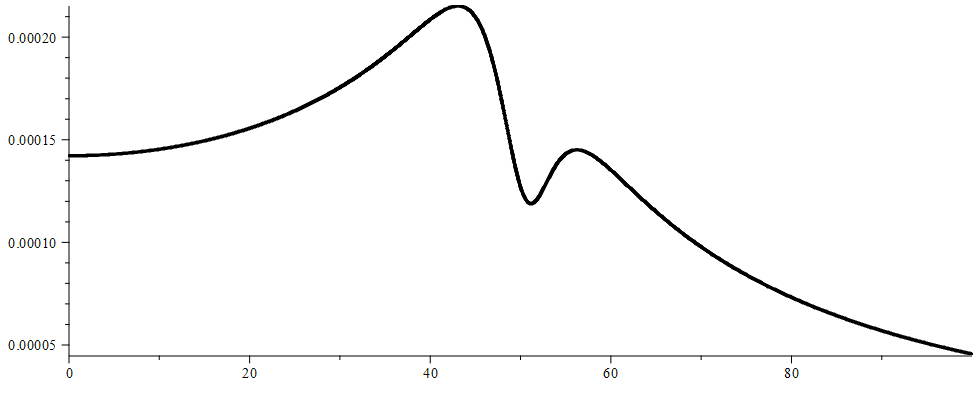}
\end{minipage}
\vspace{1em} 

\begin{minipage}{\linewidth}
  \centering
  \includegraphics[width=1\textwidth]{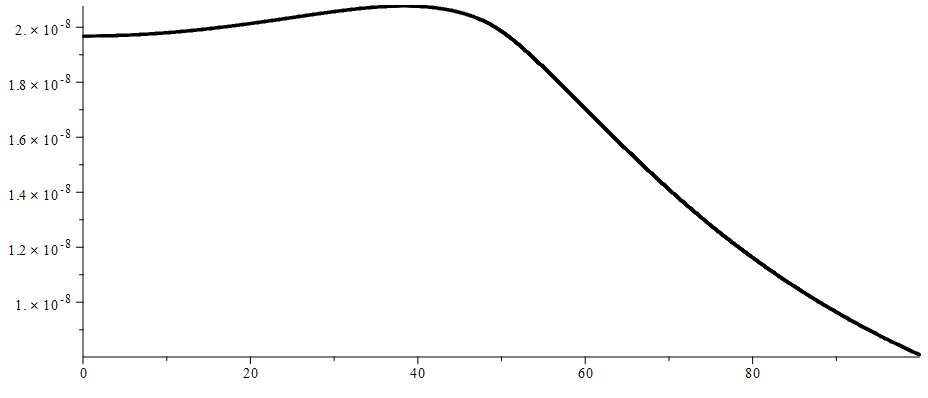}
\end{minipage}
\vspace{1em} 

\begin{minipage}{\linewidth}
  \centering
  \includegraphics[width=1\textwidth]{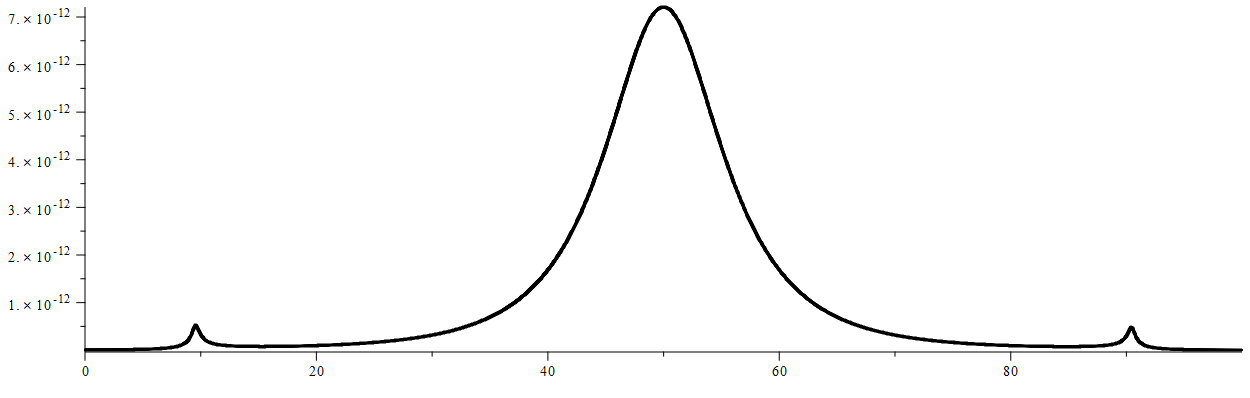}
\end{minipage}
\vspace{1em} 

\caption{Absolute difference in terms of $t$ between the approximation of the analytical Fourier transform (using splines of degree 11) and the exact analytical Fourier transform of the function $f(t)$. First figure : Boundary conditions were set $b_n=0$, $n=1,\ldots,10$, and $b_0=f(81.92)-f(0)$. Second figure : Boundary conditions were computed with Method 1. Third figure : Boundary conditions were computed analytically.}
\label{fig:figuredeg11}

\end{figure}

In Figures \ref{fig:figuredeg3} to \ref{fig:figuredeg11}, we clearly see that increasing the degree of the spline interpolation as well as improving the accuracy of the boundary conditions almost eliminates the parasitic errors and therefore significantly reduces the error on the numerical Fourier transform (\cite{beaudoin_tutorial/article_1998}). However, these figures also show that increasing the degree of the spline without improving the boundary conditions does not lead to a reduction of the error. As shown, the best scenario would be to have exact boundary conditions, but this is rarely possible in the case of a digitized function or a signal. A good compromise is to use one of the methods presented in this paper to compute the boundary conditions, which only relies on the numerical values, as they lead to very accurate results.

\section{Conclusions}
In \cite{beaudoin_tutorial/article_1998,beaudoin_new_2003}, a new spline interpolation method based on applying the DFT to a digitized function was introduced. Initially, it was believed that this method could only lead to splines of odd degrees, but further analysis showed that the method also works in the case of even degrees, as long as $N$ is chosen as an odd number. This method requires an approximation of the boundary conditions and Beaudoin and Beauchemin had suggested such an algorithm, but we have seen that it lacks robustness.

In this paper, we have proposed two new algorithms for the computation of the boundary conditions. The results were compared for splines of different degrees as well as with the traditional cubic spline. Both methods proved to be very efficient and accurate and is highly recommended when high accuracy is needed. Not only does it lead to very accurate results, but it can do so without significantly making the procedure for the computation of the spline more complex. This method can also be used to compute splines of any degree. In order to avoid the singularities that occur when $\theta$ and the number of subintervals on the mesh are even, it is however preferable to use odd-degree spline interpolation. 

Based on the results observed, we strongly believe that this method is a powerful tool that can be used to compute splines of higher degrees in order to improve the accuracy of the results when traditional methods are not sufficient. This method also provides additional information compared to the traditional methods (e.g. derivatives of any order), which is another asset of our approach. 

\section*{Acknowledgments}
We would like to acknowledge the financial support of the Natural Sciences and Engineering Research Council of Canada (NSERC) [funding reference number GC-2017-Q3-00214], the New Brunswick Innovation Foundation (NBIF), the Université de Moncton and Assumption Life. 

\bibliographystyle{plain} 
\bibliography{PepLegBeau_Splines}

\end{document}